\documentclass[12pt]{article} 
\usepackage{amssymb,latexsym}
\usepackage{amsfonts}   
\font\germ=eufm10  
\newtheorem{theorem}{Theorem}[section]
\newtheorem{pr}[theorem]{Proposition}
\newtheorem{cor}[theorem]{Corollary}
\newtheorem{de}[theorem]{Definition}
\newtheorem{con}[theorem]{Conjecture}
\newtheorem{rem}[theorem]{Remark}
\newtheorem{lem}[theorem]{Lemma}

\newtheorem{prb}[theorem]{Problem}
\newtheorem{ex}[theorem]{Example}

\newtheorem{quest}[theorem]{Question}
\def\s{\hbox{\germ S}}

\def\Z{\mathbb{Z}}

\def\Q{\mathbb{Q}}

\def\al{\alpha}
\def\d{{\Delta}}
\def\t{{\Theta}}
 
\def\q{{\cal{Q}}}

\begin{document}
\begin{center}
{\textbf{ON SOME NONCOMMUTATIVE ALGEBRAS RELATED TO K-THEORY OF FLAG 
VARIETIES, PART I}}
\end{center}

\begin{center}
\textsc{anatol n.\ kirillov and toshiaki maeno}
\end{center}
\vspace{1mm} 
\begin{center}
\small\it{ Dedicated to Adriano Garsia on the occasion of his 75-th birthday }
\vspace{5mm} 
\end{center}

\vskip -1cm
\date{May 20, 2004 }
\begin{abstract}
For any Lie algebra of classical type or type $G_2$ we define a $K$-theoretic 
analog of Dunkl's elements, the so-called truncated {\it
Ruijsenaars-Schneider-Macdonald elements}, $RSM$-elements for short, 
in the corresponding {\it Yang-Baxter 
group}, which form a commuting family of elements in the latter. For the 
root systems of type $A$ we prove that the subalgebra of the {\it bracket algebra} 
generated by the RSM-elements is isomorphic to the 
Grothendieck ring of the flag variety. In general, we prove that 
the subalgebra generated by the {\it images} of the RSM-elements 
in the corresponding {\it Nichols-Woronowicz algebra} is canonically 
isomorphic to the Grothendieck ring of the 
corresponding flag varieties of classical type or of type $G_2$. In other words,
we construct the ``Nichols-Woronowicz algebra model'' for the Grothendieck
Calculus on Weyl groups of classical type or type $G_2,$ providing a 
partial generalization of some recent results by Y.~Bazlov. We also give a 
conjectural description (theorem for type $A$) of a 
commutative subalgebra generated by the {\it truncated RSM-elements} 
in the bracket algebra for the classical root systems. Our results provide 
a proof and generalizations of recent conjecture and result by C.~Lenart 
and A.~Yong for the root system of type $A$.
\end{abstract}
\vskip 0.5cm
\nopagebreak
\begin{tabular}{rrl}
\end{tabular}


\section{Introduction}

In the paper \cite{FK} S.~Fomin and the first author have introduced a model
for the cohomology ring of flag varieties of type $A$ as a commutative 
subalgebra generated by the so-called truncated Dunkl elements in a certain
(noncommutative) quadratic algebra. This construction has been 
generalized to other root systems in \cite{KM}. The main purpose of the
present paper is to construct a $K$-theoretic analog of these 
constructions.
More specifically, we introduce certain families of pairwise commuting elements
in the Yang-Baxter group ${\cal YB}(B_n)$ or in the bracket algebra 
${\cal BE}(B_n),$ 
which conjecturally generate commutative subalgebras in the bracket algebra 
${\cal BE}(B_n)$ isomorphic to
the Grothendieck ring  of the flag varieties of type $B_n.$ 
The corresponding
results/conjectures for the flag varieties of other classical type root 
systems can be obtained
from those for the type $B$ after certain specializations. There exists the
natural surjective homomorphism
\footnote{It is believed that for a 
simply--laced (finite) Coxeter system $(W,S)$ the corresponding bracket
algebra $BE(W,S)$ and the Nichols--Woronowicz algebra ${\cal B}_{W,S}$ are
isomorphic as braided Hopf algebras. However, this is not the case in a 
non simply--laced case. For example, if $n \ge 3,$ the natural epimorphism
${\cal BE}(B_n) \longrightarrow {\cal B}_{B_n}$ has a non-trivial kernel
in degree 6. In fact, 
$Hilb({\cal BE}(B_3),t)-Hilb({\cal B}_{B_3},t)=4t^6+ \cdots.$}
from the algebra ${\cal BE}(B_n)$ to the
Nichols-Woronowicz algebra ${\cal B}_{B_n}$ of type $B.$ One of our main 
results of the present paper states that the {\it image} of our construction 
in the 
Nichols-Woronowicz algebra ${\cal B}_{B_n}$ is indeed isomorphic to the
Grothendieck ring of the flag variety of type $B_n.$ 
We also present a similar construction for the root system of type $G_2.$ 
These results can be 
viewed as a multiplicative analog/generalization for classical root 
systems and for $G_2$ of the
``Nichols-Woronowicz algebra model'' for the cohomology ring of flag varieties 
which has been constructed recently by Y.~Bazlov \cite{Ba}. 
\medskip 

In a few words the main idea behind the constructions of the paper can be
described as follows. As it was mentioned, in \cite{FK} for type $A$ and in
\cite{KM} for other root systems, a realization of the (small quantum) 
cohomology ring of flag varieties has been invented. More 
specifically, the papers mentioned above present a model for the cohomology 
ring of flag varieties as a commutative subalgebra generated 
by the so-called Dunkl elements in a certain (noncommutative) algebra. 
The main ingredient of this construction is based on
some very special solutions to the {\it classical} Yang-Baxter equation (for 
type $A$)
and {\it classical} reflection equations (for types $B,$ $C$ and $G_2$). Our 
original motivation 
was to study the related algebras and groups which correspond to the 
``quantization'' of the solutions to the classical Yang-Baxter type 
equations mentioned above, in connection with classical and quantum Schubert 
and Grothendieck Calculi. In
more detail, we define the group of ``local set-theoretical solutions'' 
to the quantum Yang-Baxter equations of type $B_n$ or of type $G_2$, together 
with the distinguished set of pairwise commuting elements in the former, the 
so-called truncated {\it Ruijsenaars-Schneider-Macdonald elements}. 
The latter is a
relativistic or multiplicative generalization of the Dunkl elements. 
For applications to the $K$-theory, we specialize the general
construction to the bracket algebra ${\cal BE}(B_n)$ and algebra 
${\cal BE}(G_2).$

Summarizing, the  main construction of our paper presents a conjectural 
description of the Grothendieck ring $K(G/B)$ corresponding to flag 
varieties $G/B$ of classical types (or $G_2$-type) to be a commutative 
subalgebra in the corresponding bracket algebra generated by the
truncated RSM-elements. To be more specific, we construct in the algebra 
${\cal BE}(B_n)$ a pairwise commuting family of elements,
{\it multiplicative or relativistic Dunkl elements}, and state a conjecture 
about the complete list of relations among the latter. 

Using some properties of the Chern homomorphism, we
prove our conjecture for the root systems of type $A.$ To our best 
knowledge,
for the root systems of type $A$ a similar description of the Grothendieck 
ring was given by C.~Lenart \cite{len1}, Lenart and Sottile \cite{len2}, 
and Lenart and Yong \cite{LY}, however without reference to 
the Yang-Baxter theory.

The main problem to prove relations between the RSM-elements in the bracket 
algebra ${\cal BE}(B_n)$ appears to be that the
intersection of kernels of all the ``braided derivations'' 
$\Delta_{ij},\Delta_{\overline{ij}}$ $1 \le i <j \le n,$ and 
$\Delta_i,$ $1\leq i \leq n,$ acting
on the algebra ${\cal BE}(B_n)$ contains only constants, see Section~5. 
At this point we 
pass to the Nichols-Woronowicz algebra ${\cal B}_{B_n}$ where the
corresponding property of the braided derivations is guaranteed, \cite{Ba}.
Since as mentioned above, there exists the natural epimorphism of braided Hopf 
algebras ${\cal B}_{B_n} \longrightarrow {\cal BE}(B_n),$  to check the 
corresponding relations in 
the Nichols-Woronowicz algebra ${\cal B}_{B_n}$ seems to be a good step to
confirm our conjectures. To prove the needed relations in the algebra
${\cal B}_{B_n},$ we develop a multiplicative analog/generalization of the
Nichols-Woronowicz algebra model for cohomology ring of flag varieties 
recently introduced by Y.~Bazlov \cite{Ba}.

In the fundamental papers \cite{KK1} and \cite{KK2} by B. Kostant and S. Kumar 
a description of the cohomology ring $H^{*}(G/B)$ and the $T$-equivariant 
$K$-theory $K_{T}(G/B)$ of a generalized flag variety $G/B$ has been obtained. 
The description of the cohomology ring $H^{*}(G/B)$ and that $K_{T}(G/B)$ 
given by Kostant and Kumar is based on the use of certain noncommutative 
algebras. The latter are suitable generalizations of the familiar nilCoxeter
$NC(W)$ and nilHecke $NH(W)$ algebras corresponding to a finite Weyl 
group $W$, to the case of Weyl groups corresponding to generalized Kac-Moody 
algebras. Note that the generators of the algebra $NC(W)$ (resp. $NH(W)$) are 
parametrized  by  the set of simple roots in the corresponding Lie algebra 
${\rm Lie}(G)$.

The main results obtained in \cite{KK1,KK2} are a far generalization of the 
well-known results in a finite dimensional case where there exists a natural 
non-degenerate  pairing between the cohomology ring $H^{*}(G/B)$ 
(resp. the Grothendieck ring $K(G/B)$) of a flag variety $G/B$ and the 
nilCoxeter algebra $NC(W)$ (resp. nilHecke algebra) of the 
Weyl group $W$ in question, see e.g. \cite{BGG}, \cite{LS1}, \cite{Ba}. In a 
few words, the pairing mentioned corresponds to a natural action of the 
Demazure operators on the cohomology ring
$H^{*}(G/B)$ (resp. the Grothendieck ring $K(G/B)$). With respect to this 
pairing the basis consisting of the Schubert polynomials in $H^{*}(G/B)$ 
(resp. the Grothendieck polynomials in $K(G/B)$) is in a duality 
with the standard basis $\{e_{w}, w \in W \},$
in the nilCoxeter algebra $NC(W)$ (resp. nilHecke algebra $NH(W)$). Under this
approach the Pieri formula for Schubert (resp. Grothendieck) polynomials is an
easy consequence of the Leibniz formula for the divided difference operators.

Our approach has its origins in the study of ``formal'' properties of the 
Pieri rules for Schubert (resp. Grothendieck) polynomials. To be more specific, the generators of our algebra corresponds to the set of positive roots in the
algebra ${\rm Lie(G)},$ and the relations are chosen 
in such a way that at first to 
guarantee the commutativity of the so-called Dunkl elements which are 
``formal'' analog of the Pieri formulas, and secondly, to guarantee the 
existence of the so-called Bruhat representation of an algebra we would like 
to construct. The existence of the Bruhat representation is a key point which 
connects our algebras with the Schubert and Grothendieck Calculi. But we 
would like to point out that our algebras have some other interesting 
representations as well, see e.g. \cite{KM}.

The above program was initiated and realized in \cite{FK} for the type $A$ 
flag varieties, and has been generalized for arbitrary finite Weyl groups in 
\cite{KM}. Another approach which is based on the theory of 
braided Hopf algebras, 
and comes up with the so-called Nichols-Woronowicz model for Schubert 
Calculus on Coxeter groups, has been developed by Y. Bazlov \cite{Ba}. One of 
the main motivations and purposes of the present paper is to construct the 
Nichols-Woronowicz model for the Grothendieck Calculus for classical Weyl 
groups and $G_2,$ as well as to generalize some  results from our previous 
paper \cite{KM} to the case of $K$-theory.

Now we want to point out the {\it main differences} between noncommutative 
algebras which have appeared in the papers by B. Kostant and S. Kumar 
\cite{KK1,KK2} and those in the present paper. First of all, the results of 
the present paper are proved only in the special case of classical root 
systems and $G_2$. The results
of \cite{KK1} and \cite{KK2} has been proved in much greater generality. 
On the other hand, in the case of type $B_n$, our algebra contains {\it 
as (dual) subalgebras} both the nilCoxeter algebra $NC(B_n)$ and a commutative 
subalgebra which is canonically isomorphic 
to the cohomology ring of $B_n$-type flag variety. Even more, 
the braided cross product of our algebra and its dual 
contains also {\it as dual subalgebras} the nilHecke algebra $NH(B_n)$ and
a commutative subalgebra which is canonically isomorphic to the Grothendieck
ring of $B_n$-type flag variety. Furthermore, an easily described 
deformation of 
our algebra contains commutative subalgebras one of which is isomorphic to 
the small quantum cohomology ring of type $B_n$-flag variety, and another one 
is conjecturally isomorphic to the quantum  $K$-theory (theorem for type $A$). 
In subsequent papers we are going to introduce (quantum) 
``degenerate affine Fomin-Kirillov'' algebra together with a commutative 
subalgebra which is isomorphic to the (quantum) $T$-equivariant $K$-theory of
type $B_n$ flag variety. We expect that our constructions can be extended to
the case of generalized flag varieties corresponding to (symmetrizable) 
Kac-Moody algebras.

Let us describe briefly the content of our paper. 

Section~2 is devoted to a general construction of commuting family of 
elements in the group ${\cal YB}(B_n)$ generated by local 
set-theoretical solutions to the family of quantum  Yang-Baxter 
equations of type $B,$ see Definition~2.1 for precise formulation. This construction
lies at the heart of our approach. In the case of type $A$ (i.e. if $g_{ij}=
h_i=1$ for all $i$ and $j$) and the Calogero-Moser representation (i.e.
$h_{ij}=1 + \partial_{ij}$) of the bracket algebra $BE(A_{n-1})$, the
elements $\Theta_1^{A_{n-1}},\cdots,\Theta_{n}^{A_{n-1}}$ correspond to the
(rational) truncated (i.e. without {\it differential} part) 
Ruijsenaars-Schneider-Macdonald operators. It seems an interesting problem
to classify all irreducible finite dimensional representations of the groups
${\cal YB}(X),$  ($X=A_{n-1},B_n,... $) together with a
simultaneous diagonalization of the operators $\Theta_{1}^{X},\cdots,
\Theta_{n}^{X}$ in these representations. 

In Section~3 we apply the result of Section~2 (Key Lemma) to construct some
distinguished  multiplicative analogue 
$\Theta_j^{A_{n-1}}:=\Theta_j^{A_{n-1}}(x)$  of the Dunkl
elements $\theta_j^{A_{n-1}},$~$1 \le j \le n,$ in the bracket algebra 
$BE(A_{n-1}).$ It happened that our elements $\Theta_j^{A_{n-1}}$ coincide 
with the $K$-theoretic Dunkl elements $1-\kappa_j$ introduced by C.~Lenart 
and A.~Yong in \cite{len1} and \cite{LY}. 
Proof of the statement that the elements
$\kappa_1,\cdots,\kappa_n$ form a family of pairwise commuting elements in
the algebra $BE(A_{n-1})$ given in \cite{len1}, appears to be quite long 
and involved. 
On the other hand, our ``Yang-Baxter approach'' 
enables us to give a simple and transparent proof that the elements
$\Theta_j^{A_{n-1}}$ mutually commute, as well as to describe relations
among these commuting elements in the algebra $BE(A_{n-1}).$ On this way we
come to the main result of Section~3, namely \medskip \\ 
{\bf Theorem A} {\it The subalgebra in $BE(A_{n-1})$ 
generated by the elements $\Theta_j^{A_{n-1}},$ $1 \le j \le n,$ is 
isomorphic to the Grothendieck ring of the flag varieties of type $A.$} 
\medskip \\
In particular, \medskip \\
{\bf Theorem B} {\it The following identity in the algebra $BE(A_{n-1})$ holds:
\[ \sum_{j=1}^{n}(\Theta_{j}^{A_{n-1}}(x))^{k}=n \] 
for any $k \in \Z.$}
\medskip \\ 
{\bf Theorem A} was stated as Conjecture~3.4 in \cite{len1} and 
Conjecture~3.5 in \cite{LY}. 
We also state a positivity conjecture as Conjecture~3.13, 
which relates the elements $\Theta_j^{A_{n-1}}$ to the Grothendieck Calculus
on the group $GL_n.$ 
This conjecture is a restatement of non-negativity conjectures from 
\cite{FK}, Conjecture~8.1, and \cite{len1}, Conjecture~3.2., 
in our setting. 

It should be emphasized that there are a lot of 
possibilities to construct a mutually commuting family of elements in the
algebra $BE(A_{n-1})$ which generate a subalgebra isomorphic to 
the Grothendieck ring $K({\cal F}l_{n}).$ For example one can take the
elements $E_1:=\exp(\theta_1^{A_{n-1}}),\cdots,E_n:= \exp(\theta_n^{A_{n-1}}).$
It is easy to see that $E_j \not= \Theta_j$ for all $j,$ however connections 
between Grothendieck polynomials and the elements $E_1,\cdots,E_n$ are not 
clear for the authors. 

Our method to describe the relations between the elements $\Theta_j^{A_{n-1}}$
is based on the study of the Chern homomorphism which relates the $K$-theory 
to the cohomology theory of flag varieties, and moreover, on description of 
the {\it commutative quotient} of the algebra $BE(A_{n-1}),
$ see Subsection~3.2. 

In Section~4 we study the $B_n$-case. First of all we introduce a modified
version ${\cal BE}(B_n)$ of the algebra $BE(B_n),$ which was introduced 
in our paper  \cite{KM}. Namely, we add additional relations in degree four, see
Definition~4.1, $(6).$  In fact we have no need to use these relations in order
to describe relations between Dunkl elements $\theta_1^{B_n},\cdots,
\theta_n^{B_n}$ in the algebra  $BE(B_n).$  However, to ensure that the $B_2$ 
Yang-Baxter relations $h_{ij}~h_{i}~g_{ij}~h_{j}=h_{j}~g_{ij}~h_{i}~h_{ij}$
are indeed satisfied, the relations $(6)$ are necessary. Another reason to 
add relations $(6)$ is that these relations are satisfied in 
the Nichols-Woronowicz 
algebra ${\cal B}_{B_{n}}.$ However, we would like to repeat again  that 
if $n \ge 3,$  then the natural 
homomorphism of algebras ${\cal BE}(B_n) \rightarrow {\cal B}_{B_{n}}$
has  a non-trivial kernel. 

The main results of Section~4 are:   

$(1)$  construction of a multiplicative 
analog $\Theta_j^{B_{n}}$ of the $B_n$-Dunkl elements 
$\theta_j^{B_{n}},$ see Definition~4.4; 

$(2)$ proof of the fact that the RSM-elements $\Theta_j^{B_{n}},$ 
$1 \le j \le n,$  form a pairwise commuting family of elements in the 
algebra  ${\cal BE}(B_n).$

Finally we give a conjectural description of all relations between the 
elements $\Theta_j^{B_{n}}.$ Here we state this conjecture in the following 
form. \medskip \\ 
{\bf Conjecture} The following identity in the algebra ${\cal BE}(B_n)$ holds:
\[ \sum_{j=1}^{n}(\Theta_j(x,y)^{B_{n}}+(\Theta_j(x,y)^{B_{n}})^{-1})^{k}=n~2^{k} \]
for all $k \in \Z_{ \ge 0}.$ \medskip 

In Section~5 we discuss on a model for the Grothendieck ring of 
flag varieties in terms of the Nichols-Woronowicz algebra for the classical 
root systems. 
Our construction is a $K$-theoretic analog of Bazlov's result \cite{Ba}. 
Our second main result proved in Section~5 is: \medskip \\ 
{\bf Theorem C} {\it Let $\varphi : B{\cal E}(B_n) \rightarrow {\cal B}_{B_n}$
be a natural homomorphism of algebras. Then
$$ \varphi(F(\Theta_1^{B_n}(x,y),\cdots,\Theta_n^{B_n}(x,y)))=0$$
for any Laurent polynomial $F$ from the defining ideal of the Grothendieck ring
of the flag variety of type $B_n.$} \medskip \\ 
{\bf Theorem C} implies the corresponding results for other classical root systems 
after some specializations. The Nichols-Woronowicz algebra ${\cal B}_{\cal X}$ 
treated in this paper is a quotient of the algebra $YB({\cal X})$ for the classical 
root system ${\cal X}.$ 
In particular, the result for $A_{n-1}$ 
is a consequence of Theorem A, 
but the argument in Section~5 is another approach based 
on the property of the Nichols-Woronowicz algebra, which works well for the 
root systems other than those of type $A_{n-1}.$ 
The idea of the proof is to construct 
the operators on the Nichols-Woronowicz algebra which induce isobaric divided 
difference operators on the commutative subalgebra generated by the RSM-elements. 

The main interest of this paper is concentrated on the classical root systems, 
for which we can use advantages of explicit handling, 
particularly in order to construct the RSM-elements. Though most of the ideas in this 
paper are expected to be applicable to an arbitrary root system, to develop 
the general framework including the exceptional root systems is a matter of concern 
for the forthcoming work. 
However, the simplest exceptional root system $G_2$ can be dealt with in similar manner 
to the case of the classical root systems. In the last section, we formulate 
the Yang-Baxter relations and define the RSM-elements for the root system of type $G_2.$ 
The argument in Section~5 again works well, so the Nichols-Woronowicz model for the 
Grothendieck ring of the flag variety of type $G_2$ is presented. 
\section{Key Lemma}
\begin{de}
Let ${\cal YB}(B_n)$ be a group generated by the elements $\{ h_{ij}, 
g_{ij} \mid 1 \le i \not= j \le n \}$ and $ \{h_i \mid 1 \le i \le n \},$
subject to the following set of relations: \smallskip 

$\bullet$ $g_{ij}=g_{ji},$ $h_{ij}=h_{ji}^{-1},$ 

$\bullet$ $h_{ij}~h_{kl}=h_{kl}~h_{ij},~~g_{ij}~g_{kl}=g_{kl}~g_{ij},$~~ 
$h_{k}~h_{ij}=h_{ij}~h_{k},$ ~~$h_{k}~g_{ij}=g_{ij}~h_{k},$ 

if all $i,j,k,l$ are distinct; \smallskip 

$h_i~h_j=h_j~h_i,$ if $1 \le i, j \le n$;
~~~$h_{ij}~g_{ij}=g_{ij}~h_{ij},$  if $1 \le i < j \le n;$ \smallskip 

$\bullet$ {\rm ($A_2$ Yang--Baxter relations)}
           $$ (1)~~~h_{ij}~h_{ik}~h_{jk}=h_{jk}~h_{ik}~h_{ij},$$
           $$ (2)~~~h_{ij}~g_{ik}~g_{jk}=g_{jk}~g_{ik}~h_{ij},$$
           $$ (3)~~~h_{ik}~g_{ij}~g_{jk}=g_{jk}~g_{ij}~h_{ik},$$
           $$ (4)~~~h_{jk}~g_{ij}~g_{ik}=g_{ik}~g_{ij}~h_{jk},$$
if $1\le i < j <k \le n;$ \smallskip 

$\bullet$  {\rm ($B_2$ quantum Yang-Baxter relation)}
$$h_{ij}~h_i~g_{ij}~h_j=h_j~g_{ij}~h_i~h_{ij},$$
if $ 1 \le i < j \le n.$
\end{de}
\begin{de}
Define the following elements in the group ${\cal YB}(B_n):$
\begin{equation}
\Theta_j=(\prod_{i=j-1}^{1}h_{ij}^{-1})~h_j~(\prod_{i=1,i \ne j}^{n}g_{ij})
~h_j~(\prod_{k=n}^{j+1} h_{jk}),
\end{equation}
for $1 \le j \le n. $ 
\end{de}
\begin{theorem}{\rm (Key Lemma )} \\

$\Theta_i \Theta_j=\Theta_j \Theta_i$ for all $ 1 \le i,j \le n.$
\end{theorem}
Our proof is based on induction plus a masterly use of the Yang-Baxter relations, 
see defining relations in the definition of the group ${\cal YB}(B_n).$
See the proof of Corollary~3.3 and Example~2.5 $(2)$ below. A complete proof 
of Theorem~2.3 one can find in Appendix. 
\begin{rem} {\rm It's not difficult to see that
$$\prod_{1 \le j \le k}\Theta_{j}=\prod_{j=1}^{k}(h_{j} \prod_{s=j+1}^{n}g_{js}
\prod_{s=1}^{j-1}g_{sj}~h_{j})~~\prod_{j=1}^{k}(\prod_{s=n}^{k+1}h_{js}).$$
In particular,
$$ \prod_{j=1}^{n} \Theta_{j}= (\prod_{k=1}^{n} 
(\prod_{j \le k}g_{jk})~h_{k})^2 .$$ }
\end{rem}
\begin{ex} {\rm $(1)$  Take $n=2$. Then $\Theta_1=h_1~g_{12}~h_1~h_{12}$ and
$\Theta_2=h_{12}^{-1}~h_2~g_{12}~h_2.$ 
Let us check that $\Theta_1$ and $\Theta_2$ commute. Indeed, using the 
$B_2$-quantum Yang-Baxter relation $h_{12}~h_1~g_{12}~h_2=h_2~g_{12}~h_1~h_{12}$ and 
the commutativity relation $h_1~h_2=h_2~h_1,$ we see that
$$\Theta_1~\Theta_2=h_1~g_{12}~h_1~h_2~g_{12}~h_2=
h_{12}^{-1}~(h_{12}~h_1~g_{12}~h_2)h_1~g_{12}~h_2$$
$$=h_{12}^{-1}~h_2~g_{12}~h_1(h_{12}~h_1~g_{12}~h_2)=h_{12}^{-1}~h_2~g_{12}~h_1~h_2~g_{12}~h_1~h_{12}=
\Theta_2~\Theta_1=(h_1~g_{12}~h_2)^2.$$ 
$(2)$ Take $n=3$. Then we have 
$$\Theta_1=h_1~g_{12}~g_{13}~h_1~h_{13}~h_{12}, \; 
\Theta_2=h_{12}^{-1}~h_2~g_{12}~g_{23}~h_2~h_{23}, \; 
\Theta_3=h_{23}^{-1}~h_{13}^{-1}~h_3~g_{13}~g_{23}~h_3,$$ 
and 
$$\Theta_1~\Theta_2~\Theta_3=(h_1~g_{12}~h_2~g_{13}~g_{23}~h_3)^2.$$
Let us illustrate the main ideas behind the proof of Key Lemma by the 
following example. \medskip \\ 
$ \Theta_{1}~\Theta_{3}~\Theta_{1}^{-1}=h_1~g_{12}~g_{13}~h_1~h_{13}~
{\bf h_{12}~h_{23}^{-1}~h_{13}^{-1}}~h_3~{\bf g_{13}~g_{23}}~h_3~
{\bf h_{12}^{-1}}~h_{13}^{-1}~h_1^{-1}~g_{13}^{-1}~g_{12}^{-1}~h_1^{-1}$ \smallskip 

$ = h_1~{\bf g_{12}~g_{13}}~h_1~{\bf h_{23}^{-1}}~h_3~
g_{23}~{\bf g_{13}~h_3~h_{13}^{-1}~h_1^{-1}}g_{13}^{-1}~g_{12}^{-1}~h_1^{-1}$ ~(by (1) 
and (2)) \smallskip 

$= h_1~h_{23}^{-1}~g_{13}~{\bf g_{12}}~h_3~
{\bf g_{23}~h_{13}^{-1}}~h_3~g_{12}^{-1}~h_1^{-1}$ ~(by (4) and $B_2$-YBE) \smallskip 

$ = h_1~h_{23}^{-1}~{\bf g_{13}~h_3~h_{13}^{-1}}~g_{23}~
h_3~{\bf h_1^{-1}}$ ~(by (3)) \smallskip 

$= \Theta_{3}$ ~(by $B_2$-YBE). }
\end{ex}
\begin{rm} {\rm  We define the groups ${\cal YB}(A_{n-1})$ and 
${\cal YB}(D_n)$ to be the quotients of that ${\cal YB}(B_n)$ by
the normal subgroups generated respectively by the elements  $\{h_i,g_{ij},
1 \le i < j \le n \}$ and $\{h_i, 1 \le i \le n \}.$  The group 
${\cal YB}(G_2)$ will be defined in Section~6. We expect that the
subgroup in ${\cal YB}(B_n)$ generated by the elements $\Theta_1^{B_n},
\cdots,\Theta_{n}^{B_n}$ is isomorphic to the free abelian group of rank $n.$
It seems an interesting problem to construct analogs of the group
${\cal YB}(B_n)$ and the elements $\Theta_1^{B_n},
\cdots,\Theta_{n}^{B_n}$ for any (finite) Coxeter group.
}
\end{rm}
\begin{quest} {\rm Does there exist a finite-dimensional {\it faithful} 
representation of the group ${\cal YB}(X),$ $X = A_{n-1},B_n,...$ ?} 
\end{quest}

\section{Algebras  $YB(A_{n-1})$ and $BE(A_{n-1})$}

\subsection{Definitions and main results}
(i) {\bf Algebra $YB(A_{n-1})$}

\begin{de} Let $R$ be a $\Q$-algebra. 
Define the algebra $YB_R(A_{n-1})$ as an associative algebra over $R$ 
generated by the elements $h_{ij}(x),$ $1 \le i \ne j \le n,$ 
$x\in R,$ subject to the relations $(0) - (4):$ \medskip 

$(0)$  $h_{ij}(x)h_{ji}(x)=1,$ \smallskip

$(1)$  $h_{ij}(x)h_{ij}(y)=h_{ij}(x+y);$ in particular, 
$h_{ij}(x)h_{ij}(-x)=1,$ \smallskip 

$(2)$  $h_{ij}(x)h_{kl}(y)=h_{kl}(y)h_{ij}(x),$ if $i,j,k,l$ are distinct, \smallskip

$(3)$ $h_{ij}(x)h_{jk}(y)+h_{ik}(x+y)=h_{jk}(y)h_{ik}(x)+
h_{ik}(y)h_{ij}(x),$ \smallskip

$h_{jk}(y)h_{ij}(x)+h_{ik}(x+y)=h_{ik}(x)h_{jk}(y)+h_{ij}(x)h_{ik}(y),$

if $1 \le i < j < k \le n,$ \smallskip

$(4)$ $h_{ik}(x)~(h_{ij}(x)-h_{ik}(y))~h_{ij}(y)=h_{ij}(y)~(h_{ij}(x)-h_{ik}(y))~h_{ik}(x),$

if $1 \le i < j < k \le n.$

\end{de}
For any element $z\in R$ we denote by  $YB(A_{n-1})[z]$ (resp. $YB(A_{n-1})$) the algebra over $\Q$ generated
by the elements $h_{ij}(z)$ and $h_{ij}(-z),$ (resp. $h_{ij}(1)$ and $h_{ij}(-1)$),
$1 \le i \ne j \le n.$

\begin{lem}{\rm (Quantum Yang-Baxter equation)} \\
The following relations in the algebra 
$YB(A_{n-1})[z]$
\begin{equation}
h_{ab}(z)h_{ac}(z)h_{bc}(z)=h_{bc}(z)h_{ac}(z)h_{ab}(z),
~~1 \le a < b < c \le n,
\end{equation}
in the algebra $YB(A_{n-1})[z]$ 
are a consequence of the relations $(0)-(4).$
\end{lem}
\begin{cor} Define elements  $\Theta_j^{A_{n-1}}(z),$ $j=1,\ldots, n,$ in 
the algebra $YB(A_{n-1})[z]$ as follows:
\begin{equation}
\Theta_{j}^{A_{n-1}}(z)=h_{j-1,j}^{-1}(z) \cdots h_{1j}^{-1}(z)~h_{jn}(z) 
\cdots h_{j,j+1}(z), ~~1 \le j \le n.
\end{equation}
Then 
\[ \Theta_{j}^{A_{n-1}}(z)\Theta_{k}^{A_{n-1}}(z)=\Theta_{k}^{A_{n-1}}(z)
\Theta_{j}^{A_{n-1}}(z), \; \; \; {\rm for} \; \; {\rm all} \; \; 1 \le j,k \le n. \]
\end{cor}
This Corollary is a particular case of Key Lemma above. We would like to 
include a separate proof of this special case to show the main ideas behind 
the usage of the Yang-Baxter relations, and since in this case the proof 
is much easy.

{\it Proof.} It is enough to check that if $1 \le i \le j \le n,$ then 
$$ \Theta_{i} ~\Theta_{j} ~\Theta_{i}^{-1}=\Theta_j.$$
By definition, \smallskip 

$\Theta_{i} ~\Theta_{j} ~\Theta_{i}^{-1}=h_{i-1,i}^{-1} \cdots 
h_{1,i}^{-1}~h_{i,n} \cdots h_{i,i+1}~h_{j-1,j}^{-1}\cdots {\bf h_{i+1,j}^{-1}}~{\bf h_{i,j}^{-1}}~h_{i-1,j}^{-1} \cdots h_{1,j}^{-1}$ 
\smallskip \\ 
\hspace{30mm}
$h_{j,n} \cdots h_{j,j+1}~
{\bf h_{i,i+1}^{-1}}h_{i,i+2}^{-1}\cdots h_{i,n}^{-1}~h_{1,i} \cdots h_{i-1,i}.$ 
\smallskip \\ 
Using local commutativity relations, see Definition 3.1 $(2),$ we can move the
factor ${\bf h_{i,i+1}^{-1}}$ to the left till we have touched on the factor 
${\bf h_{i,j}^{-1}}$. As a result, we will come up with the triple product:
$${\bf h_{i+1,j}^{-1}~h_{i,j}^{-1}~h_{i,i+1}^{-1}},$$
which is equal, according to the Yang-Baxter relation $(3.2),$ to the product
$${\bf h_{i,i+1}^{-1}~h_{i,j}^{-1}~h_{i+1,j}^{-1}}.$$
Now we can move the factor ${\bf h_{i,i+1}^{-1}}$ to the left to cancel it
with the term  $h_{i,i+1},$ which comes from the rightmost factor in 
the element $\Theta_i.$

As a result, we will have 
\[ \Theta_{i} ~\Theta_{j} ~\Theta_{i}^{-1}= 
h_{i-1,i}^{-1}\cdots h_{i,i+2}~h_{j-1,j}^{-1}\cdots {\bf h_{i+2,j}^{-1}}~{\bf h_{i,j}^{-1}} \cdots h_{j,j+1}^{-1}~{\bf h_{i,i+2}^{-1}} \cdots h_{i-1,i}. \]
Now we can move to the left the factor ${\bf h_{i,i+2}^{-1}}$ till we have
touched on the factor ${\bf h_{i,j}^{-1}}$ to give the triple product
$$ {\bf h_{i+2,j}^{-1}~h_{i,j}^{-1}~h_{i,i+2}^{-1}},$$ 
which is equal to ${\bf h_{i,i+2}^{-1}~h_{i,j}^{-1}~h_{i+2,j}^{-1}}.$ Now we can move the factor ${\bf h_{i,i+2}^{-1}}$ to the left to cancel 
it with the corresponding factor $h_{i,i+2},$ and so on.

It is readily  seen that finally we will come to the element $\Theta_j.$  
\rule{3mm}{3mm}  \\ 
It is clear that $\prod_{j=1}^{n} \Theta_{j}^{A_{n-1}}(z)=1.$ 
\begin{rem} {\rm Let $\Theta_j(z):=\Theta_{j}^{A_{n-1}}(z).$ Then it is not true 
that $\Theta_j(x)\Theta_k(y)= \Theta_k(y)\Theta_j(x),$ if 
$j \ne k,$ $x \ne y.$ }
\end{rem}
\begin{rem}
{\rm Though the algebra $YB(A_{n-1})[z]$ can be constructed as 
a quotient of the group algebra $\Q \langle {\cal YB}(A_{n-1}) \rangle ,$ 
they are not isomorphic.}
\end{rem}
\begin{theorem}{\rm (Main theorem, the case of algebra $YB(A_{n-1})[z]$)}
\begin{equation}
\prod_{j=1}^{n}(1+(1-\Theta_{j}^{A_{n-1}}(z))t)=1.~~~Equivalently, 
~~\prod_{j=1}^{n}(1+\Theta_{j}^{A_{n-1}}(z)t)=(1+t)^n.
\end{equation}
\end{theorem}
This theorem is equivalent to:
\begin{theorem} Let  $G_j^{A_{n-1}}=\Theta_j^{A_{n-1}}(z)-1,$ 
$1 \le j \le n.$   Then, after the substitution $z=1,$
$$ e_{j}(G_{1}^{A_{n-1}},\ldots,G_{n}^{A_{n-1}})=0, ~~~ 1 \le j \le n $$ 
is the complete list of relations in the algebra $YB(A_{n+1})$ 
among the elements $G_1^{A_{n-1}},\ldots,G_n^{A_{n-1}}.$ 
Here, $e_j$ is the $j$-th elementary symmetric polynomial. 
\end{theorem}
The proof is given in Subsection 3.2. 
It is based on the properties of the {\it Chern homomorphism}. 
  
\begin{cor} The algebra over $\Z$ generated by the elements $G_{1}^{A_{n-1}}
\vert_{z=1},\ldots,G_{n}^{A_{n-1}} \vert_{z=1}, $ is canonically isomorphic to 
the {\it integral} Grothendieck ring $K({\cal F}l_n)$  of the flag manifold 
of type $A_{n-1}.$ 
\end{cor} 
(ii) {\bf Algebra $BE(A_{n-1})$ }
\begin{de}{\rm (\cite{FK})}  Define algebra $BE(A_{n-1})$ 
$($denoted by ${\cal E}_n$ in 
{\rm \cite{FK}}$)$ as an associative  algebra over 
$\Z$ with generators $x_{ij},$ $1 \le i \ne j \le n,$ subject to the following relations 
\smallskip 

$(0)$ $x_{ij}+x_{ji}=0,$  $1 \le i \ne j \le n,$

$(1)$ $x_{ij}^{2}=0,$  $1 \le i \ne j \le n,$

$(2)$ $x_{ij}~x_{jk}+x_{jk}~x_{ki}+x_{ki}~x_{ij}=0,$ if all $i,j,k$ are 
distinct.
\end{de}

The Dunkl elements $\theta_j,$ $j=1,\ldots,n,$ in the algebra 
$BE(A_{n-1})$ are defined by $ \theta_j:=\theta_j^{A_{n-1}}=\sum _{i \ne j}x_{ij}.$

The Dunkl elements 
form a pairwise commuting family of elements in the algebra $BE(A_{n-1}),$ 
\cite{FK},  and generate a commutative subalgebra
in $BE(A_{n-1}),$ which is canonically isomorphic to the cohomology ring
$H^{*}({\cal F}l_n)$ of the flag variety ${\cal F}l_n$ of type $A_{n-1},$ 
\cite{FK}. 

For an element $t$ of a $\Q$-algebra $R,$ 
define $h_{ij}(t)=1+tx_{ij}=\exp(tx_{ij}) \in BE(A_{n-1}) \otimes R.$ 
\begin{lem} The elements $h_{ij}(t), 1 \le i,j \le n,$ satisfy the all relations
$(0)-(4)$ of the definition of the algebra $YB(A_{n-1}).$
\end{lem}
We will use the same notation $\Theta_j^{A_{n-1}}$, ~$1 \le j \le n,$  to 
denote the elements in the algebra 
$BE(A_{n-1})$ defined by the formula (3.3). It follows from Corollary 3.3 that 
they form a pairwise commuting family of elements in the algebra $BE(A_{n-1}).$

It's clear that
$\Theta_{j}^{A_{n-1}}(z)=1+z~\theta_{j}^{A_{n-1}}+\cdots$, and the product
in the RHS of (3.3) may  be written as follows:
\begin{equation}
\Theta_j^{A_{n-1}}(z)= \sum (-1)^s x_{b_{1},j}~x_{b_{2},j} \cdots x_{b_{s},j}~
x_{j,a_1}~x_{j,a_2} \cdots x_{j,a_r}~z^{r+s},
\end{equation}
where the sum runs over the all sequences of integers $(a_1 > a_2 > \cdots > 
a_r)$ and $(b_1 > b_2 > \cdots > b_s)$ such that $n \ge a_1 > a_r > j > b_1 >
b_s \ge 1;$ cf. \cite[Section 2]{len1}.

Remember that $G_j^{A_{n-1}}:=\Theta_j^{A_{n-1}}-1,$~$1 \le j \le n.$
\begin{de} Let $w \in S_n$ be a permutation. 
Define the Grothendieck polynomial
${\cal G}_{w}(X_n) \in \Z~[X_n]$ to be a unique polynomial of the form
${\cal G}_{w}(X_n)=\sum_{\al \subset \delta_n}c_{\al}(w)~x^{\al}$ such that
\begin{equation}
{\cal G}_{w}(G_1^{A_{n-1}},\ldots,G_n^{A_{n-1}}) \cdot {\it id} = w 
\end{equation}
in the Bruhat representation  of the algebra $BE(A_{n-1})$ 
{\rm (}see  {\rm \cite[Section 3.1]{FK}} 
{\rm )}, where $\delta_n : =(n-1,n-2,\ldots,1,0)$ and $X_n:=(x_1,\ldots,x_n).$
\end{de}
It is not difficult to see that the Grothendieck polynomials 
defined here coincide with those introduced in \cite{L1}, see also 
\cite{len2}. 
\begin{cor} Let $u \in S_n$ and $v \in S_n$ be two permutations. Assume that
in the group ring $\Z \langle S_n \rangle$ of the symmetric group $S_n$ 
one has the following equality: 
$$ {\cal G}_{u}(G_1^{A_{n-1}},\ldots,G_n^{A_{n-1}})\cdot v = \sum_{w \in S_n}
c_{u,v}^{w}~w. $$
Then the coefficient $c_{u,v}^{w}$ is equal to the multiplicity of the
Grothendieck polynomial ${\cal G}_{w}(X_n)$ in the product of ${\cal G}_{u}(X_n)$
and ${\cal G}_{v}(X_n):$
$${\cal G}_{u}(X_n)~{\cal G}_{v}(X_n)=\sum_{w \in S_n} c_{u,v}^{w}~{\cal G}_{w}(X_n)$$
in the Grothendieck ring $K({\cal F}l_n)$ of the flag manifold of type 
$A_{n-1}.$
\end{cor}

\begin{con} {\rm For any permutation $w \in S_{n}$ the value of the Grothendieck
polynomial ${\cal G}_{w}(x_1,\ldots,x_n)$ after the substitution  
$x_1:=G_{1}^{A_{n-1}}, \ldots, x_n:=G_{n}^{A_{n-1}},$ and $z=1$, 
can be written as a linear combination of 
monomials in  $x_{ij}$'s , $1 \le i < j \le n,$ with {\bf non-negative}
integer coefficients.}
\end{con}
\begin{ex}{\rm (Grothendieck-Pieri formula in the algebra $BE(A_{n-1})$, cf
\cite{len2}) 
$$1+{\cal G}_{(k,k+1)}(G_1,\ldots,G_n)=\prod_{1 \le j \le k}\Theta_j =
\prod_{j=1}^{k}\prod_{s=n}^{k+1}h_{js}=\sum \prod_{j=1}^{r}x_{a_{j},b_{j}},$$
where the sum runs over all sequences of integers $(1 \le a_1 \le \cdots \le 
a_r \le k)$ and $(b_1,\ldots,b_r)$ such that $k < b_j \le n$, $j=1,\ldots, r,$
~and $a_i=a_{i+1} \Rightarrow b_i > b_{i+1}.$}
\end{ex}
Our methods allow to obtain a {\bf subtraction free} formula in the algebra 
$BE(A_{n-1})$ for the value of the Grothendieck polynomials 
${\cal G}_{(k,k+1,\cdots,k+r)}(G_1,\ldots,G_n),$ ~$1 \le k \le n-r-1,$ as well. ~We hope to report on our results in a separate publication. 
\begin{ex} {\rm Take $n=3,$ then 

$\Theta_1:=\Theta_{1}^{A_2}(1)=h_{13}(1)~h_{12}(1)=
1+x_{12}+x_{13}+x_{12}~x_{13},$ 

$\Theta_2:=\Theta_{2}^{A_2}(1)= 
h_{12}^{-1}(1)~h_{23}(1)=1-x_{13}+x_{23}-x_{13}~x_{12}-x_{23}~x_{13},$

$\Theta_3:=\Theta_{3}^{A_2}(1)=h_{23}^{-1}(1)~h_{13}^{-1}(1)=1-x_{13}-x_{23}+
x_{23}~x_{13}.$ \smallskip \\
As a preliminary step, we compute the elementary symmetric polynomials 
$e_{k}(\Theta_1,
\Theta_2,\Theta_3),$ $k=1,2,3.$ Indeed, it's easily seen from the formulae
above that 
$\Theta_1+\Theta_2+\Theta_3=3$ and $\Theta_1~\Theta_2~\Theta_3=1.$
To compute $e_2(\Theta_1,\Theta_2,\Theta_3),$ all one has to do is to apply
the following relation
$$h_{12}~h_{23}^{-1}=h_{23}^{-1}~h_{13}+h_{13}^{-1}~h_{12}-1,$$
where we put by definition $h_{ij}:=h_{ij}(1).$ The former equality follows
from the relation (3) in Definition 3.1. 
Hence, 
\[ e_2(\Theta_1,\Theta_2,\Theta_3)=h_{13}~h_{23}+h_{13}~{\bf h_{12}~h_{23}^{-1}}~h_{13}^{-1}+h_{12}^{-1}~h_{13}^{-1} \] 
\[ =h_{13}~h_{23}+h_{13}~h_{23}^{-1}+
h_{12}~h_{13}^{-1}-1+h_{12}^{-1}~h_{13}^{-1}=2h_{13}+2h_{13}^{-1}-1=3. \]
To continue, let us list the Grothendieck polynomials ${\cal G}_{w}(x)$ 
corresponding to the symmetric group $S_3$:
\[ {\cal G}_{id}(x)=1, \; {\cal G}_{s_1}(x)=x_1, \; {\cal G}_{s_2}(x)=
x_1+x_2+x_1x_2, \] 
\[ {\cal G}_{s_1s_2}(x)=x_1x_2, \; {\cal G}_{s_2s_1}(x)=
x_{1}^{2}, \; {\cal G}_{w_{0}}(x)=x_{1}^{2}x_{2}. \]
Now let us consider the substitution $x_{j}=G_{j}=\Theta_{j}(1)-1,$ 
$j=1,2,3.$ More explicitly, $G_1=x_{12}+x_{13}+x_{13}~x_{12}$ and $G_2=-x_{12}+
x_{23}-x_{13}~x_{12}-x_{23}~x_{13}.$ Therefore,
\[ {\cal G}_{s_2}(G_1,G_2)=x_{13}+x_{23}+x_{13}~x_{23}, \; {\cal G}_{s_1s_2}
(G_1,G_2)=x_{13}~x_{23}+x_{23}~x_{13}, \] 
\[ {\cal G}_{s_2s_1}(G_1,G_2)= x_{12}~x_{13}+x_{13}~x_{12}, \] 
\[ {\cal G}_{w_{0}}(G_1,G_2)=x_{12}~x_{13}~x_{23}+x_{13}~x_{12}~x_{13}+x_{13}~x_{23}~x_{13}+x_{13}~x_{12}~x_{13}~x_{23}. \]
Finally, let us consider the commutative subalgebra in $BE(A_2) \otimes \Q$ 
generated by the  elements $E_j:=\exp(\theta_j),$ $j=1,2,3.$ It's not 
difficult to check that
$$ 2E_1= h_{13}~h_{12}+h_{12}~h_{13},~2E_2=h_{12}^{-1}~h_{23}+h_{23}~h_{12}^{-1},~2E_3=h_{23}^{-1}~h_{13}^{-1}+h_{13}^{-1}~h_{23}^{-1}.$$
It is an easy matter as well to see that the subalgebra in $BE(A_2) \otimes 
\Q$ ~generated over $\Q$ by 
the elements  $E_i,$ $i=1,2,3,$ is isomorphic to the algebra 
$\Q[\Theta_1,\Theta_2,\Theta_3].$ 
In particular, for all symmetric polynomials $f(x_1,x_2,x_3)$ we have
$$ f(1-E_1,1-E_2,1-E_3)=0.$$}
\end{ex}

\begin{pr} The subalgebra in $BE(A_{n-1}) \otimes \Q$ generated by the 
elements  $E_i:=\exp(\theta_i),$ $1 \le i \le n,$ is isomorphic to the algebra
over $\Q$ generated by the elements $\Theta_j^{A_{n-1}},$ $1 \le j \le n.$

In particular, the complete list of relations among the elements $1-E_1,
\ldots,1-E_n$ in the quadratic algebra $BE(A_{n-1})$ is given by
$$ e_i(1-E_1,\ldots,1-E_n)=0,$$
for $i=1,\cdots,\ n.$ Thus the commutative subalgebra generated by the 
elements
$\exp(\theta_1),\ldots,\exp(\theta_n)$ is isomorphic 
to the rational Grothendieck ring $K({\cal F}l_n) \otimes \Q$ of the flag 
manifold ${\cal F}l_n$ of type $A_{n-1}.$ 
\end{pr}
However, it seems that there are no direct connections of the elements $E_j$'s
with the Grothendieck Calculus.

\begin{rem}{\rm More generally, let $Q(t) \ne 0$ be a polynomial such that $Q(0)=0.$
Define the elements $q_{i}:=1+Q(\theta_i),$  $1 \le i \le n,$ in the algebra

$BE(A_{n-1}).$ It's clear that the elements $q_1,\ldots,q_n$ pairwise
commute, and 
$$ e_i(q_{1}-1,\ldots,q_{n}-1)=0, ~~1 \le i \le n.$$
}
\end{rem}

\begin{rem}{\rm (Quantum Grothendieck Calculus)} \\ 
{\rm It is easy to see that the relations in Definition 3.1 are still true, if
we replace the condition $(1)$ in Definition 3.9 by the following one

$(1')$ $x_{ij}^2=q_{ij}$, $1 \le i < j \le n,$ where the parameters $q_{ij}$
are assumed to commute with all the generators $x_{kl}$, 
$1 \le k < l \le n.$

The algebra over $\Z \lbrack\ q_{ij} \mid 1 \le i < j \le n  \rbrack $ generated by the elements  $x_{ij},$ $ 1 \le i \ne j \le n,$ 
subject to the
relations $(0)$, $(1')$ and $(2),$ is called the {\it quantized bracket 
algebra} and denoted by $qBE(A_{n-1}),$ cf. \cite[Section 15]{FK} and \cite{Kir2}.   

As a corollary we see that the elements $\Theta_j^{q},$  $1 \le j \le n,
$ defined by the formula $(3.2),$ form a pairwise commuting family of elements
in the algebra $qBE(A_{n-1}).$ 

After the specialization 
\[ q_{ij}= \left\{ 
\begin{array}{cc}
q_i, & \textrm{if $j=i+1,$} \\ 
0, & \textrm{otherwise,} 
\end{array}
\right. \]
the multiplicative Dunkl elements generate the quantum Grothendieck ring 
in the sense of Givental and Lee \cite{GL}. The generalization to the equivariant $K$-theory 
is an open problem.}
\end{rem}

\begin{prb} {\rm 
Describe the commutative subalgebras in the quantized algebra 
$qYB(A_{n-1})$ generated by

$(1)$ $ \Theta_1^{q}(1), \ldots,  \Theta_n^{q}(1),$

$(2)$ ${\widetilde E}_1:=\exp(\theta_1), \ldots, {\widetilde E}_n:=
\exp(\theta_n).$}
\end{prb}

\subsection{Chern homomorphism}
Denote by ${\cal H}:=BE(A_{n-1})^{{\it ab}}\otimes \Q$ the quotient of the algebra
$BE(A_{n-1})$ by its commutant. 
It is known, \cite[Proposition 4.2]{FK}, that the algebra 
$BE(A_{n-1})^{{\it ab}}$ has dimension $n!$, and its Hilbert polynomial is
given by
$$ Hilb(BE(A_{n-1})^{{\it ab}},t)=(1+t)(1+2t) \cdots (1+(n-1)t).$$
Denote by $1+{\cal H}^{+}$ the
multiplicative monoid generated by the elements of the form $1+h$, where 
$h \in {\cal H}$ does not have the term of degree zero. 
\begin{pr} Let $R^{(n-1)}$ be the subspace of 
the commutative subalgebra \\ 
$R=\Q[\theta_1,\ldots,\theta_n] \subset 
BE(A_{n-1})\otimes \Q$ whose elements are of degree $\leq n-1.$ Then the 
subspace $R^{(n-1)}$ is injectively mapped into ${\cal H}$ by 
the quotient homomorphism $BE(A_{n-1})\otimes \Q \rightarrow {\cal H}.$ 
\end{pr}
{\it Proof.} Since the algebra $R$ is isomorphic to the coinvariant algebra of 
the symmetric group, the monomials 
\[ \theta_1^{i_1}\cdots \theta_{n-1}^{i_{n-1}}, \; \; \; 0\leq i_k \leq n-k, \]
form a linear basis of $R.$ 
The linear map $R^{(n-1)}\rightarrow {\cal H}$ induced by the quotient homomorphism 
is a homomorphism between $S_n$-modules. Hence, it is enough to show the images 
of the monomials $\theta_1^{i_1}\cdots \theta_{n-1}^{i_{n-1}}$ do not vanish 
in ${\cal H}$ for $(i_1,\ldots,i_{n-1})$ such that $\sum_{k=1}^{n-1}i_k=n-1$ and 
$i_1\geq i_2 \geq \cdots \geq i_{n-1}.$ 
We expand the monomials $\theta_1^{i_1}\cdots \theta_{n-1}^{i_{n-1}}$ 
of this form in the algebra $BE(A_{n-1})\otimes \Q$ 
by using the Pieri formula proved by Postnikov \cite{Po}, 
(first conjectured in \cite{FK}). 
The Pieri formula shows that 
\[ e_k(\theta_1,\ldots,\theta_m)=\widetilde{\sum}[i_1j_1]\cdots[i_kj_k], \] 
where $\widetilde{\sum}$ stands for the multiplicity-free sum, and 
$(i_1,j_1),\ldots,(i_k,j_k)$ run over all pairs such that $i_a \leq m < j_a \leq n,$ 
$a=1,\ldots,k,$ and all $i_a$'s are distinct. 

On the other hand, the monomials of form 
\[ [i_1j_1]\cdots[i_kj_k], \; \; \; i_a < j_a \; \; (a=1,\ldots,k), \; \; 
j_1<j_2<\cdots < j_k, \]
give a linear basis of ${\cal H}$ (\cite[Corollary 10.3]{Kir2}). By the involution 
$\omega : [ij] \mapsto [n+1-j \; , \; n+1-i],$ we have a linear basis of form 
\begin{equation} 
[i_1j_1]\cdots[i_kj_k], \; \; \; i_a < j_a \; \; (a=1,\ldots,k), \; \; 
i_1<i_2<\cdots < i_k. 
\end{equation} 
For each monomial expression $[i_1j_1]\cdots[i_kj_k]$ in ${\cal H},$ we define 
\[ \mu ([i_1j_1]\cdots[i_kj_k]) := \sum_{m=1}^k (j_m-i_m). \] 
Every element in ${\cal H}$ can be expressed as a linear combination of the monomials 
listed in (3.7) by repeatedly applying the substitution 
$[ab][ac] \rightarrow [ab][bc] - [ac][bc]$ with $a<b<c.$ 
On each step of the procedure, the monomials of minimal $\mu$ appearing in the expression 
of $\theta_1^{i_1}\cdots \theta_{n-1}^{i_{n-1}}$ with $i_1+\cdots +i_{n-1}=n-1$ are 
not cancelled or are replaced by new ones. 
So one can check the image of $\theta_1^{i_1}\cdots \theta_{n-1}^{i_{n-1}}$ in 
${\cal H}$ is not zero. 
\rule{3mm}{3mm} 
\begin{de} Define the Chern homomorphism $($to the commutative quotient$)$ 
$$ c': YB(A_{n-1}) \rightarrow 1+{\cal H}^{+} $$
by the following rules: \\

$\bullet$  $c'(f+g)=c'(f)c'(g),$ ~if ~$f,g \in YB(A_{n-1}),$

$\bullet$  $c'(\prod_{i < j} h_{ij}^{n_{ij}})=1+\sum_{i < j}n_{ij}~x_{ij}.$
\end{de}
It is clear that $c'(\Theta_j)=1+\theta_j,$ $\forall j.$ 
\begin{rem} {\rm 
We can also define the homomorphism 
\[ c: \Q[\Theta_1,\ldots,\Theta_n] \rightarrow \Q[\theta_1,\ldots,\theta_n] \] 
by the conditions $c(f+g)=c(f)c(g)$ and $c(\Theta_j)=1+\theta_j,$ $j=1,\ldots,n,$ 
which is compatible with 
the Chern homomorphism (in the usual sense) 
\[ c: K(Fl_n) \rightarrow 1+H^+(Fl_n). \] 
However, the homomorphism $c'$ defined above does not coincide with $c$ in the part 
of degree $\geq n.$ Indeed, the maximal degree of the commutative quotient 
${\cal H}$ is $n-1.$} 
\end{rem} 
\begin{pr} {\rm (cf. \cite[Section 5]{LS1})}
For any permutation $w \in S_n,$
$$c(1+{\cal G}_{w}(G_1,\ldots,G_n))=1-(-1)^{l(w)}(l(w)-1)!~{\s}_{w}
(\theta_1,\ldots,\theta_n)+\sum_{u}a_{u}(w)~{\s}_{u}(\theta_1,\ldots,\theta_n),$$
where the sum ranges over all permutations $u \in S_n$ such that $l(u) > l(w),$ 
and $a_u(w)$ is a constant in ${\bf Z}$ determined by $u$ and $w.$ 
\end{pr}
{\it Proof of Theorem 3.7.} Note that the commutative quotient $YB(A_{n-1})^{ab}$ 
is isomorphic to the algebra ${\cal H}.$ Moreover, The subspace of polynomials 
of degree $\leq n-1$ in the 
RSM-elements $\Theta_1^{A_{n-1}},\ldots ,\Theta_n^{A_{n-1}}$ in $YB(A_{n-1})$ 
is also injectively mapped into ${\cal H}$ from Proposition 3.20. 
We regard $1+{\cal H}^+$ as an $(n!-1)$-dimensional $\Q$-linear 
space so that the homomorphism 
$\bar{c}:{\cal H}^+ \rightarrow 1+{\cal H}^+$ induced by 
the Chern homomorphism $c'$ is a $\Q$-linear map. The image of the 
linear basis (3.7) of ${\cal H}^+$ by the 
homomorphism $\bar{c}$ is linearly independent. Hence, 
$\bar{c}:{\cal H}^+\rightarrow 1+{\cal H}^+$ is an isomorphism 
between linear spaces. Since it is easy to see 
\[ c'(e_j(\Theta_1^{A_{n-1}},\ldots , \Theta_n^{A_{n-1}}))=1 \in 1+{\cal H}^+, 
\; \; \; 1\leq j \leq n-1, \] 
one can conclude that 
\[ e_{j}(G_{1}^{A_{n-1}},\ldots,G_{n}^{A_{n-1}})=0, ~~~ 1 \le j \le n-1. \] 
The equality 
\[ \prod_{i=1}^n \Theta_i^{A_{n-1}}=1 \] 
in the algebra $YB(A_{n-1})$ can be obtained by direct computation. 
\rule{3mm}{3mm} 
\begin{prb} {\rm 
Construct a lift of $c'$ to 
\[ YB(A_{n-1}) \rightarrow 1+BE(A_{n-1})^+ \] 
in some suitable sense.} 
\end{prb}
\section {Algebras ${\cal BE}(B_n)$ and $YB(B_n)$}
(i) {\bf Algebra ${\cal BE}(B_n)$} (cf. \cite{KM})
\begin{de} Define the algebra ${\cal BE}(B_n)$ as the algebra 
$($say, over $\Q)$  with
generators
$$ [i,j],~\overline{[i,j]},~ 1 \le i \ne j \le n, ~~{\rm and} ~~[i], ~1 \le i \le n, 
$$
subject to the following relations: 
\smallskip \\ 
$(0)$ $[i,j]=-[j,i],$ $\overline{[i,j]}=\overline{[j,i]},$ \\ 
$(1)$ $[i,j]^2=0$, $\overline{[i,j]}^2=0,$ $1 \le i < j \le n,$ and $[i]^2=0$,
$1 \le i \le n,$ 
\smallskip\\
$(2)$ $[i,j][k,l]=[k,l][i,j],$ 
$\overline{[i,j]}[k,l]=[k,l]\overline{[i,j]},$ 
$\overline{[i,j]}\overline{[k,l]}=\overline{[k,l]}
\overline{[i,j]},$ \\ 
\hspace*{5.2mm} if $\{ i,j\} \cap \{ k,l\} =\emptyset$, 
\smallskip \\ 
$(3)$ $[i][j]=[j][i],$ $[i,j]\overline{[i,j]}=
\overline{[i,j]}[i,j],$ 
$[i,j][k]=[k][i,j]$, if $k\not= i,j$, 
\smallskip \\ 
$(4)$\hspace*{5.2mm}  $[i,j][j,k]+[j,k][k,i]+[k,i][i,j]=0,$ \\ 

\hspace*{5.2mm} $\overline{[i,k]}[i,j]+[j,i]\overline{[j,k]}
+\overline{[k,j]}\overline{[i,k]}=0,$ \\ 

\hspace*{5.2mm} $[i,j][i]+[j][j,i]+[i]\overline{[i,j]}+
\overline{[i,j]}[j]=0,$ \\
\hspace*{5.2mm} if all $i,$ $j$ and $k$ are distinct,
\smallskip \\  
$(5)$ $[i,j][i]\overline{[i,j]}[i]+\overline{[i,j]}[i][i,j][i]+
[i][i,j][i]\overline{[i,j]}+[i]\overline{[i,j]}[i][i,j]=0,$  if $i<j,$ 
\bigskip \\
$(6)$ $[i,j][i]\overline{[i,j]}[j]=[j]\overline{[i,j]}[i][i,j],$ if $i<j.$
\end{de}
\begin{rem} {\rm (a) In the definition of the algebra $BE(B_n),$ see \cite[Section 9.1]{KM}, 
the condition $(6)$ is absent. In fact, there is no need to use 
the latter condition for the purposes of \cite{KM}. However, we need the 
condition $(6)$ to ensure the $B_2$ quantum Yang-Baxter relation, 
which is necessary for our construction of  a commutative family of elements 
in the algebra  $YB(B_n),$ see (ii) below. \smallskip \\ 
(b) In \cite{KM}, the authors has introduced the quantum deformation 
$qBE(B_n)$ of the bracket algebra. Similarly we introduce the quantum 
deformation of the algebra $q{\cal BE}(B_n)$  
which is generated by the same symbols as in 
${\cal BE}(B_n)$ and is obtained by replacing 
the relation in $(1)$ corresponding 
to the simple roots by \smallskip 

$[i,i+1]^2=q_i$, $1 \le i \le n-1,$ and $[n]^2=q_n.$ \smallskip \\
In the subsequent construction, we can work in the quantum bracket 
algebra $q{\cal BE}(B_n)$ instead of ${\cal BE}(B_n).$ The RSM-elements 
in Definition 4.4 also form a commuting family of elements in 
$q{\cal BE}(B_n).$ 
Though it is expected that the RSM-elements in the quantum setting should 
describe the quantum Grothendieck Calculus in $B_n$-case, 
the relations satisfied by them in the algebra $q{\cal BE}(B_n)$ 
are not clearly seen.} 
\end{rem}
The Dunkl elements  \cite{KM} are given by 
\begin{equation} {\theta}_i := \theta^{B_n}_i
=\sum_{j\not= i}([i,j]+\overline{[i,j]})+2[i], \; \; \;  1\leq i \leq n.
\end{equation}
Note that  the Dunkl elements 
$\tilde{\theta}_i$ correspond to the Monk type formula in the 
cohomology ring of the flag variety of type $B.$
\bigskip \\ 
(ii) {\bf Algebra $YB(B_n)$}
\medskip \\ 
Let $x$ and $y$ be elements in a $\Q$-algebra $R.$ 
Define the algebra $YB(B_n)$ as a subalgebra in ${\cal BE}(B_n)\otimes R$ 
generated over $R$ by the elements: 

$h_{ij}:=\exp(x[i,j])=1+x[i,j], ~~g_{ij}:=
\exp(x \overline{[i,j]})=1+x \overline{[i,j]},$ $1 \le i < j \le n,$ 

and $h_j:= \exp(y[j])=1+y[j],$ $1 \le j \le n.$ 
\begin{pr} The elements $h_{ij}, g_{ij} ~and ~h_k, ~1 \le i <j \le n, ~1 \le k \le n,$ 
satisfy the all relations listed in Definition 2.1.
\end{pr}

\begin{de}  Define
$$\Theta_{j}^{B_n}(x,y)=(\prod_{i=j-1}^{1}h_{ij}(x)^{-1})~h_j(y)~(\prod_{i=1,i \ne j}^{n}g_{ij}(x))~h_j(y)~(\prod_{k=n}^{j+1} h_{jk}(x)),$$
for $1 \le j \le n. $ 
\end{de}
\begin{cor} The elements $\Theta_{j}^{B_n}(x,y)$ commute pairwise.
\end{cor} 
\begin{rem} {\rm  It is not difficult to see that 
$$\Theta_{j}^{B_n}(1,1) \not= \exp (\theta_j^{B_n}), $$
where $\theta_j^{B_n}, 1 \le j \le n,$ denote the $B_n$-Dunkl elements in the
algebra ${\cal BE}(B_n)$. The commuting family of elements 
$\exp(\theta_j^{B_n}), 1 \le j \le n,$ also generate a (finite dimensional)
commutative subalgebra in ${\cal BE}(B_n) \otimes \Q$. However, we don't 
know the complete list of relations among these elements. }
\end{rem}
\begin{con}{\rm (The case of algebra  $YB(B_n)$)} \\
In the algebra $YB(B_n)$ we have the
following identity
\begin{equation}
\prod_{j=1}^{n}(1+(\Theta_{j}^{B_{n}}(x,y)+
(\Theta_{j}^{B_{n}}(x,y))^{-1})t)= (1+2t)^{n}.
\end{equation}
Equivalently,
\begin{equation}
\prod_{j=1}^{n}(1+\Theta_{j}^{B_n}(x,y)t)(1+(\Theta_{j}^{B_n}(x,y))^{-1}t)=
(1+t)^{2n}.
\end{equation}
\end{con}
This conjecture is equivalent to: 
\begin{con} Let $G_{j,\alpha}^{B_{n}}=(\Theta_{j}^{B_{n}}(x,y))^{\alpha}-(\Theta_{j}^{B_{n}}(x,y))^{- \alpha},$ 
$1 \le j \le n,$ ~$\alpha \in \Q.$
Then
\begin{equation}
e_j((G_{1,\alpha}^{B_{n}})^2,\ldots,(G_{n,\alpha}^{B_{n}})^2)= 0, ~1 \le j \le n.
\end{equation}
\end{con}
\begin{rem} {\rm Theorem 3.6, i.e. the equality 
\begin{equation}
\prod_{j=1}^{n}(1+\Theta_{j}^{A_{n-1}}~t)=(1+t)^{n}.
\end{equation}
follows from Conjecture 4.7.}
\end{rem}
{\it Proof.} The multiplicative Dunkl elements $\Theta_{j}^{A_{n-1}}(x)$ can 
be obtained from those  $\Theta_{j}^{B_{n}}(x,y)$ after the specialization 
$y:=0$ and $g_{ij}:=1.$ Since $\prod_{j=1}^{n} \Theta_{j}^{A_{n-1}}=1$, 
it follows from Conjecture 4.7 that if we denote by
$P_{n}(t)$ the LHS of (4.11) then $P_{n}(t)P_{n}(t^{-1})=
(1+t)^{n}(1+t^{-1})^{n}.$ Therefore, $P_{n}(t)=(1+t)^{n}.$  \rule{3mm}{3mm}


\begin{rem}
{\rm The algebra $YB(C_n)$ can be naturally identified with the algebra $YB(B_n).$ 
The corresponding RSM-elements relate via 
\[ \Theta_{j}^{C_n}(x,y)=\Theta_{j}^{B_n}(x,y/2). \]
}
\end{rem}

\section{Nichols-Woronowicz model for Grothendieck ring of flag varieties}
In the preceeding sections, we have tried to construct the models of the 
Grothendieck ring $K(G/B)$ in the algebras $BE(A_{n-1})$ and 
${\cal BE}(B_n)$ for the corresponding root systems respectively. 
The algebras $BE(A_{n-1})$ and ${\cal BE}(B_n)$ have braided Hopf algebra 
structures. In particular, $BE(A_{n-1})$ is conjecturally isomorphic 
to the so-called Nichols-Woronowicz algebra. Bazlov \cite{Ba} has 
constructed the model of the coinvariant algebra of the finite Coxeter 
group from this viewpoint. In this section, we construct a model 
of the Grothendieck ring of the flag variety in terms of the 
Nichols-Woronowicz algebra associated to a Yetter-Drinfeld module over 
the Weyl group $W$ for the classical root systems. 

The Nichols-Woronowicz algebra ${\cal B}(V)$ is a braided Hopf algebra 
determined by a given braided vector space $V=(V,\psi).$ The braided 
vector space $(V,\psi)$ is a finite-dimensional vector space 
$V$ equipped with the braiding $\psi: V \otimes V 
\rightarrow V\otimes V$ that is a canonically given linear endomorphism 
satisfying the braid relation 
\[ \psi_{12}\psi_{23}\psi_{12}= \psi_{23}\psi_{12}\psi_{23}: 
V\otimes V \otimes V \rightarrow V\otimes V \otimes V, \] 
where $\psi_{ij}:V\otimes V \otimes V \rightarrow V\otimes V \otimes V$ 
is obtained by applying $\psi$ on the $i$-th and $j$-th components. 
The Nichols-Woronowicz algebra ${\cal B}(V)$ is a braided 
analog of the symmetric tensor algebra, which is defined by replacing 
the symmetrizer by the braided symmetrizer. 
For an element $w\in S_n$ with a reduced decomposition 
$w=s_{i_1}\cdots s_{i_l},$ the linear endomorphism $\psi_w:=
\psi_{i_1 \, i_1+1}\cdots \psi_{i_l\, i_l+1}$ on $V^{\otimes n}$ is 
well-defined from the braid relation. The Woronowicz symmetrizer 
$\sigma_n(\psi): V^{\otimes n} \rightarrow V^{\otimes n}$ is given by 
the formula 
\[ \sigma_n(\psi):= \sum_{w\in S_n} \psi_w. \]
The Nichols-Woronowicz algebra ${\cal B}(V)$ is the quotient of 
the tensor algebra $\oplus_n V^{\otimes n}$ by the kernels of 
the braided symmetrizers $\sigma_n(\psi)$: 
\[ {\cal B}(V)= \mathop{\oplus}_n V^{\otimes n}/ \mathop{\oplus}_n {\rm Ker}
(\sigma_n(\psi)) . \] 
The Nichols-Woronowicz algebra ${\cal B}(V)$ provides a natural framework 
to perform the braided differential calculus. 

Let us consider the Nichols-Woronowicz algebra ${\cal B}_{\cal X}$ obtained from 
the Yetter-Drinfeld module 
\[ V=\mathop{\oplus}_{\alpha\in \Psi}\Q[\alpha]/
([\alpha]+[-\alpha])_{\alpha\in \Psi} \]
for 
the root system $\Psi$ of classical type ${\cal X},$ 
(${\cal X}=A_{n-1},$ $B_n,$ $C_n,$ $D_n$). 
Let $W({\cal X})$ be the corresponding Weyl group. 
The $W({\cal X})$-action on $V$ is given by $w([\alpha])=[w(\alpha)],$ 
and the $W({\cal X})$-degree of $[\alpha]$ is the reflection $s_{\alpha}.$ 
The structure of the braided vector space on $V$ is given by the braiding 
$\psi([\alpha]\otimes [\beta])=[s_{\alpha}(\beta)]\otimes [\alpha].$ 
For the details on the definition of the algebra ${\cal B}_{\cal X},$ see \cite{Ba}. 
The algebra ${\cal B}_{\cal X}$ is a quotient of the algebra $YB({\cal X}).$ 

The Weyl group $W(B_n)$ acts on the algebra $YB(B_n).$ Denote by 
$s_1=s_{12},\ldots ,s_{n-1}=s_{n-1 \, n},$ and $s_n$ the simple 
reflections. The subgroup $S_n=W(A_{n-1})\subset W(B_n)$ acts 
on $YB(B_n)$ via the permutation of the indices of $h_{ij},$ 
$g_{ij}$ and $h_i.$ The action of the simple reflection $s_n$ 
is given as follows: 
\[ s_n(h_{ij})= h_{ij}, \; s_n(g_{ij})=g_{ij}, \; s_n (h_i)= h_i, 
\; \; \; {\rm for} \; \; i,j\not= n, \]
\[ s_n (h_{in})=g_{in}, \; s_n (g_{in})=h_{in}, \; s_n (h_n)= h_n^{-1}. \]
Define the twisted derivations $\d_{ij}$ $(i<j)$ and $\d_i$ on 
$YB(B_n)$ by 
\[ \d_{ij}(h_{kl})= \left\{
\begin{array}{cc} 
1, & {\rm if} \; \; i=k \; \; {\rm and} \; \; j=l, \\
0, & {\rm otherwise,} 
\end{array}
\right. \]
\[ \d_{ij}(g_{ij}) = \d_{ij}(h_k)=0, \] 
\[ \d_{i}(h_j)= \left\{
\begin{array}{cc} 
1, & {\rm if} \; \; i=j, \\
0, & {\rm otherwise,} 
\end{array}
\right. \]
\[ \d_{i}(h_{jk}) = \d_{i}(g_{jk})=0, \] 
and the twisted Leibniz rule 
\[ \d_{ij}(xy)= \d_{ij}(x)y+s_{ij}(x)\d_{ij}(y), \] 
\[ \d_i(xy) = \d_i (x)y+s_i(x)\d_i(y). \] 
Let us consider the operators $\q_i := h_{i \, i+1}^{-1} \circ \d_{i \, i+1}$ 
$(i<n)$ and $\q_n := h_n^{-1} \circ \d_n$ on $YB(B_n).$ 
\begin{lem}
Let $\Theta_{j}:=\Theta_{j}^{B_{n}}(1,1).$ One has 
\[ \q_i(\t_j)= \left\{ 
\begin{array}{cc} 
\t_{i+1}, & {\rm if} \; \; j=i, \\ -\t_{i+1} & {\rm if} \; \; j=i+1, \\
0, & {\rm otherwise,} 
\end{array} \right. \] 
for $i<n,$ and 
\[ \q_n (\t_j) = \left\{ 
\begin{array}{cc} 
1+\t_n^{-1}, & {\rm if} \; \; j=n, \\
0, & {\rm otherwise.} 
\end{array} \right. \] 
\end{lem}
{\it Proof.} It is clear that $\q_i(\t_j)=h_{i \, i+1}^{-1}\d_{i\, i+1}(\t_j)=0$ 
$(i<n)$ for 
$j\not= i,i+1$ and $\q_n (\t_j) = h_n^{-1} \d_n (\t_j)=0$ for $j\not=n.$ We have 
by direct computation 
\[ \q_i (\t_i)= h_{i \, i+1}^{-1}\d_{i \, i+1}\left( \prod_{k=i-1}^1 h_{ki}^{-1}\cdot h_i
\prod_{k=1,k\not= i}^n g_{ki}\cdot h_i \cdot \prod_{k=n}^{i+1} h_{ik} \right) \]
\[ = h_{i \, i+1}^{-1}\cdot \left( \prod_{k=i-1}^1 h_{k\, i+1}^{-1}\cdot h_{i+1}
\prod_{k=1,k\not= i+1}^n g_{k\, i+1}\cdot h_{i+1} \cdot \prod_{k=n}^{i+2} h_{i+1\, k} 
\right) = \t_{i+1} , \] 
\[ \q_i (\t_{i+1})= h_{i \, i+1}^{-1}\d_i \left( \prod_{k=i}^1 h_{k\, i+1}^{-1}
\cdot h_{i+1} \prod_{k=1,k\not= i+1}^n g_{k\, i+1}\cdot h_{i+1} \cdot 
\prod_{k=n}^{i+2} h_{i+1\, k} \right) \]
\[ = h_{i \, i+1}^{-1}\cdot \left( -\prod_{k=i-1}^1 h_{k\, i+1}^{-1}\cdot h_{i+1}
\prod_{k=1,k\not= i+1}^n g_{k\, i+1}\cdot h_{i+1} \cdot \prod_{k=n}^{i+2} h_{i+1\, k} 
\right) = -\t_{i+1} . \] 
Similarly, 
\[ \q_n (\t_n)= h_n^{-1} \d_n \left( \prod_{k=n-1}^1 h_{kn}^{-1} \cdot h_n 
\prod_{k=1}^{n-1} g_{kn} \cdot h_n \right) \]
\[ = h_n^{-1} \left( h_n +  \prod_{k=n-1}^1 g_{kn}^{-1} \cdot h_n^{-1} 
\prod_{k=1}^{n-1} h_{kn} \right) = 1+ \t_n^{-1}.  \; \; \; \rule{3mm}{3mm} \] 
\begin{lem} The simple reflections act on the elements $\t_1,\ldots, \t_n$ 
as follows. 
\[ h_{i\, i+1}^{-1}\cdot s_i(\t_j) \cdot h_{i\, i+1} =
\t_{s_i(j)}, \; \; \; {\rm for} \; \; i=1,\ldots, n-1, \]
\[ h_n^{-1}\cdot s_n(\t_j) \cdot h_n= \left\{ \begin{array}{cc}
\t_n^{-1}, & {\rm if} \; \; j=n, \\
\t_j, & {\rm otherwise.} 
\end{array} \right. \]
\end{lem}
{\it Proof.} 
If $j\not=i,i+1,$ then the equality 
\[ h_{i\, i+1}^{-1}\cdot s_i(\t_j) \cdot h_{i\, i+1} =
\t_j \] 
follows from the relations 
\[ h_{i\, i+1}^{-1}h_{ji}h_{j\, i+1}h_{i\, i+1}=h_{j\, i+1}h_{ji} \] 
and 
\[ h_{i\, i+1}^{-1}g_{i+1\, j}g_{ji}h_{i\, i+1}=g_{ij}g_{i+1\, j}. \] 
We also have 
\begin{eqnarray*}
\lefteqn{h_{i\, i+1}^{-1}\cdot s_i(\t_i) \cdot h_{i\, i+1}} \\
& = & h_{i\, i+1}^{-1} \left( \prod_{k=i-1}^1h_{k\, i+1}^{-1}\cdot h_{i+1}
\prod_{k=1,k\not=i+1}^n g_{k\, i+1} \cdot h_{i+1} \cdot \prod_{k=n}^{i+2}
h_{i+1\, k}\cdot h_{i+1\, i} \right)\cdot h_{i\, i+1} \\
& = & \t_{i+1} , 
\end{eqnarray*} 
and this completes the proof of the first equality. 

We can obtain 
\[ h_n^{-1}h_{jn}h_jg_{jn}h_n=g_{jn}h_jh_{jn} \] 
from the $B_2$ Yang-Baxter relation. 
This shows the equality 
\[ h_n^{-1}\cdot s_n(\t_j) \cdot h_n = \t_j \] 
for $j\not=n.$ Since 
\[ s_n(\t_n)= \prod_{k=n-1}^1g_{kn}^{-1}\cdot h_n^{-1}\prod_{k=1}^{n-1}h_{kn}\cdot 
h_n^{-1}, \] 
we have 
\[ h_n^{-1}\cdot s_n(\t_n) \cdot h_n = \t_n^{-1}. \; \; \; \; \rule{3mm}{3mm} \]

Consider the action of $W(B_n)$ on the ring of Laurent polynomials 
$\Q[X_1^{\pm 1},\ldots,X_n^{\pm 1}]$ via 
\[ (w f)(X_1,\ldots,X_n):=f(X_{w(1)},\ldots,X_{w(n)}), \; \; \; w\in S_n=W(A_{n-1}), \] 
and 
\[ (s_n f)(X_1,\ldots,X_n):=f(X_1,\ldots,X_{n-1},X_n^{-1}). \] 
\begin{lem}
Let $F(\t)$ and $G(\t)$ be Laurent polynomials in $\t_1,\ldots,\t_n.$ Then, 
\[ \q_i (F(\t)G(\t))= \q_i(F(\t))G(\t)+ (s_iF)(\t)\q_i(G(\t)), \; \; \; i=1,\ldots,n . \]
\end{lem}
{\it Proof.} The equalities in Lemma 5.2 imply 
\[ h_{i\, i+1}^{-1}\cdot s_i (F(\t))\cdot h_{i\, i+1} = (s_iF)(\t), \] 
so 
\begin{eqnarray*}
\q_i (F(\t)G(\t)) & = & h_{i\, i+1}^{-1} \d_{i\, i+1} (F(\t)G(\t)) \\ 
& = & h_{i\, i+1}^{-1} \d_{i\, i+1} (F(\t))G(\t)+h_{i\, i+1}^{-1}s_i(F(\t))h_{i\, i+1}
\cdot h_{i\, i+1}^{-1}\d_{i\, i+1}(G(\t)) \\
& = & \q_i(F(\t))G(\t)+ (s_iF)(\t)\q_i(G(\t)) 
\end{eqnarray*} 
for $i<n.$ The equality 
\[ \q_n (F(\t)G(\t))= \q_n(F(\t))G(\t)+ (s_nF)(\t)\q_n(G(\t)) \] 
is proved in the same way. \rule{3mm}{3mm} 

Define the operators $\tau_1,\ldots,\tau_{n-1}$ and $\tau_{n}:=
\tau_{n}^{B_{n}}$ on  $\Q[X_1^{\pm 1},\ldots,X_n^{\pm 1}]$ 
by 
\[ (\tau_i f)(X):= X_{i+1}\frac{f(X)-(s_if)(X)}{X_i-X_{i+1}}, \; \; \; i=1,\ldots,n-1, \] 
\[ (\tau_n f)(X):= \frac{f(X)-(s_nf)(X)}{X_n-1}. \] 
The operator corresponding to $\tau_n$ in the case of type $C_n$ is given by 
\[ (\tau_{n}^{C_n}f)(X):= \frac{f(X)-(s_nf)(X)}{X_n^2-1}. \]
We consider the group $W(D_{n})$ as the subgroup of $W(B_{n}).$ 
Let $\tau_{n}^{D_n}:=
\tau_{n}^{B_n}\tau_{n-1}\tau_{n}^{B_n}.$ Then we have 
\[ (\tau_{n}^{D_n}f)(X_1,\cdots,X_{n-1},X_n)= \frac{f(X_1,\cdots,X_{n-1},X_n)-
f(X_1,\cdots,X_n^{-1},X_{n-1}^{-1})}{X_{n-1}X_{n}-1}.
\]
\begin{pr} Let $\t_j:=\Theta_{j}^{B_n}(1,1),$ ~$1 \le j \le n,$ then
\[ \q_i(F(\t_1,\ldots,\t_n))=(\tau_i F)(\t_1,\ldots,\t_n) . \]
\end{pr}
{\it Proof.} This follows from Lemmas 5.1 and 5.3. \rule{3mm}{3mm}
\medskip \\
\begin{rem} 
{\rm One can obtain the corresponding results for $A_{n-1}$ (resp. $D_n$) after 
specialization $g_{ij}=h_i=1$ (resp. $h_i=1$), $\forall i,j.$ }
\end{rem}
\begin{rem}
{\rm All the construction in this section till Proposition 5.4 can be 
done on the level of the group algebra $\Q \langle {\cal YB}(B_n) \rangle.$} 
\end{rem}
We have the homomorphisms 
\[ \varphi : YB(A_{n-1}) \rightarrow BE(A_{n-1}) \rightarrow {\cal B}_{A_{n-1}}, \] 
\[ \varphi : YB(D_n) \rightarrow BE(D_n) \rightarrow {\cal B}_{D_n}, \] 
\[ \varphi : YB(B_n) \rightarrow BE(B_n) \rightarrow {\cal B}_{B_n}, \] 
given by $h_{ij} \mapsto 1+[ij],$ $g_{ij} \mapsto 1+\overline{[ij]}$ and 
$h_i \mapsto 1+[i].$ \\
Conjecturally, the quadratic algebras $BE(A_{n-1})$ and $BE(D_n)$
\cite{KM} are isomorphic respectively to the Nichols-Woronowicz
algebras ${\cal B}_{A_{n-1}}$ and ${\cal B}_{D_n}.$  

The Nichols-Woronowicz algebra is equipped with the duality pairing 
\[ \langle \; , \; \rangle : {\cal B}_{\cal X} \otimes {\cal B}_{\cal X} \rightarrow \Q \] 
and naturally defined braided derivations acting on it. 
Here we are interested in the derivations 
$\overline{D}_{[\alpha]}$ given by the formula 
\[ \overline{D}_{[\alpha]}(\xi)=({\rm id}_{\cal B} \otimes \langle \; , \; \rangle)
(\psi_{V,{\cal B}} \otimes {\rm id}_{\cal B})([\alpha]\otimes \xi_{(1)} \otimes \xi_{(2)}), \]
where $\psi_{V,{\cal B}}: V \otimes {\cal B} \rightarrow {\cal B}\otimes V$ is the 
braiding induced by $\psi,$ and we use Sweedler's notation $\triangle (\xi)= 
\xi_{(1)} \otimes \xi_{(2)}$ for the coproduct $\triangle$ of the Nichols-Woronowicz 
algebra. 
The twisted derivations $\d_{ij},$ $\d_{\overline{ij}}$ and $\d_i$ are 
corresponding to the derivations on the Nichols algebras, 
namely $\varphi (\d_{ij}(x))=\overline{D}_{ij}(\varphi(x)),$ 
$\varphi (\d_{\overline{ij}}(x))=\overline{D}_{\overline{ij}}(\varphi(x)),$ 
$\varphi (\d_i(x))=\overline{D}_i(\varphi(x)).$ 
Note that the intersection of the kernels of all the derivations 
$\overline{D}_{[\alpha]}$ coincides with the degree zero part ${\cal B}^0_{\cal X}=\Q.$ 
This is the essential property of the Nichols-Woronowicz algebra which will be 
used in the subsequent argument. 

Let $P$ be the weight lattice associated to some root system and 
$\Q[P]=\Q[e^{\lambda} | \lambda \in P]$ 
its group algebra. Denote by $\epsilon: \Q[P] \rightarrow \Q$ 
the algebra homomorphism given by $e^{\lambda}\mapsto 1,$ $\forall \lambda\in P.$ 
The Grothendieck ring of the corresponding flag variety can be 
expressed as a quotient algebra $\Q[P]/I,$ where the ideal $I$ is generated by 
the $W$-invariant elements of form $f-\epsilon (f).$ 
\begin{theorem}
Let $F$ be a Laurent polynomial in the defining ideal of the Grothendieck ring 
of the flag variety of classical type ${\cal X},$ and $\t_j:= \Theta_{j}^{\cal X}.$ Then, 
\[ \varphi(F(\t_1,\ldots,\t_n))=0  \] 
in the corresponding Nichols-Woronowicz  algebra 
${\cal B}_{\cal X}$ $({\cal X}=A_n,$ $B_n,$ $C_n$ or $D_n).$
\end{theorem}
{\it Proof.} In the following, we consider the root system of type $B_n.$ The
cases of type $A,C,D$ can be obtained from this case by a certain specialization.
For simplicity, we use the same symbol $\t_i$ for 
the corresponding element to the RSM-elements in ${\cal B}_{B_n}.$ 
Let $\epsilon_j(X):=e_j(X_1+X_1^{-1},\ldots,X_n+X_n^{-1})-e_j(2,\ldots,2).$ 
Proposition 5.4 implies that 
\[ \varphi( \q_i(\epsilon_j(\t)))=0. \] 
Hence, we have $\overline{D}_{i\, i+1}(\epsilon_j(\t))=0$ and 
$\overline{D}_n (\epsilon_j(\t))=0.$ From the $W$-invariance of 
the polynomial $\epsilon_j$ and Lemma 5.2, 
it follows that $s_k(\epsilon_j(\t))=
h_{k\, k+1}\epsilon_j(\t)h_{k\, k+1}^{-1}$ and 
$s_n(\epsilon_j(\t))=
h_n \epsilon_j(\t)h_n^{-1}.$ 
Thus, for $k\not= i,$ 
\[ \d_{i\, i+1}(s_k(\epsilon_j(\t)))=\d_{i\, i+1}(h_{k\, k+1}\epsilon_j(\t)h_{k\, k+1}^{-1}) \] 
\[ =s_i(h_{k\, k+1})\d_{i\, i+1}(\epsilon_j(\t))h_{k\, k+1}^{-1}=0. \] 
For $k=i,$ 
\[ \d_{i\, i+1}(s_i(\epsilon_j(\t)))=\d_{i\, i+1}(h_{i\, i+1}\epsilon_j(\t)h_{i\, i+1}^{-1}) \]
\[ =(\epsilon_j(\t)- h_{i\, i+1}^{-1}\cdot s_i(\epsilon_j(\t))\cdot h_{i\, i+1})h_{i\, i+1}^{-1} + 
h_{i\, i+1}^{-1}\d_{i\, i+1}(\epsilon_j(\t))h_{i\, i+1}^{-1}=0. \]
More generally, one can show that if $\d_{kl}(\epsilon_j(\t))=0,$ then 
$\d_{kl}(s_i(\epsilon_j(\t)))=0.$ Since $w\circ \d_{kl} \circ w^{-1} = \d_{w(k) \, w(l)}$ 
for $w\in W,$ we can conclude that 
$\overline{D}_{kl}(\epsilon_j(\t))=\overline{D}_{\overline{kl}}(\epsilon_j(\t))=0,$ $\forall k,l.$ 
Similarly, $\overline{D}_k(\epsilon_j(\t))=0,$ $\forall k.$ 
Since the constant term of $\epsilon_j(\t)$ considered as a polynomial in $[ab]$'s, 
$\overline{[ab]}$'s and $[a]$'s is zero, 
it follows that $\epsilon_j(\t)=0$ in ${\cal B}_{B_n}.$ \rule{3mm}{3mm} 

Let us remark that it follows from the above considerations that 
in the case of $D_n$ we have relations 
$e_{k}(\Theta_{1}^{D_n}+(\Theta_{1}^{D_n})^{-1}, \cdots,\Theta_{n}^{D_n}+(\Theta_{n}^{D_n})^{-1})=0$ for $1 \le k < n$ and the additional relation
$\prod_{j=1}^{n}((\Theta_{j}^{D_n})^{1/2}-(\Theta_{j}^{D_n})^{-1/2})=0$ in ${\cal B}_{D_n}.$ 

For any Laurent polynomial $F$ that is not in the ideal generated by $\epsilon_1,\ldots, 
\epsilon_n,$ one can find a sequence of indices $i_1,\ldots , i_r$ such that 
$\tau_{i_1}\cdots \tau_{i_r}F(X)\in \Q \setminus \{ 0 \}.$ Hence we have the following. 
\begin{cor}
The RSM elements $\Theta_i$ generate the algebra isomorphic to the Grothendieck ring 
of the corresponding flag variety of classical type as a commutative subalgebra 
in ${\cal B}_{\cal X}.$ 
\end{cor}
\begin{rem}
{\rm The operators $\q_1,\ldots , \q_n$ satisfy the relations 
\[ \q_i^2=\q_i, \; i=1,\ldots n, \] 
\[ \q_i \q_{i+1} \q_i = \q_{i+1} \q_i \q_{i+1}, \; i=1,\ldots , n-2, \] 
\[ \q_{n-1} \q_n \q_{n-1} \q_n = \q_n \q_{n-1} \q_n \q_{n-1}. \] 
}
\end{rem}
\section{The case of root system of type $G_2$}
Let us consider the root system of type $G_2.$ Let 
\[ \Psi_+= \{ a,\; b, \; c, \; d, \; e, \; f \}  \] 
be the set of positive roots, 
where $a$ and $f$ are the simple roots and $b=3a+f,$ $c=2a+f,$ $d=3a+2f,$ 
$e=a+f.$ 
\begin{de}Denote by ${\cal YB}(G_2)$ the group generated by six 
elements $h_a,h_b,h_c,h_d,h_e,h_f$ subject to the following relations:

$\bullet$  $h_a~h_d=h_d~h_a,~~h_b~h_e=h_e~h_b,~~h_c~h_f=h_f~h_c;$

$\bullet$ $(A_2$ - Yang-Baxter relation$)$ ~~~~$h_b~h_d~h_f=h_f~h_d~h_b;$

$\bullet$ $(G_2$ - Yang-Baxter relation$)$
$$ h_a~h_b~h_c~h_d~h_e~h_f=h_f~h_e~h_d~h_c~h_b~h_a.$$
\end{de}
\begin{pr}Define the $RSM$-elements of type $G_2$ in  ${\cal YB}(G_2)$
as follows
$$\Theta_1^{G_2}:=h_d~h_b~h_c~h_d~h_e~h_f, \; \; \; \Theta_2^{G_2}:=
h_f^{-1}~h_b~h_d~h_c~h_b~h_a.$$
Then we have $\Theta_1^{G_2}~\Theta_2^{G_2}=\Theta_2^{G_2}~\Theta_1^{G_2}.$
\end{pr} 
Let us consider the group algebra $\Q \langle {\cal YB}(G_2) \rangle .$ 
The Weyl group $W(G_2)$ naturally acts on the algebra 
$\Q \langle {\cal YB}(G_2) \rangle .$ 
The twisted derivations $\Delta_a$ and $\Delta_f$ determined by the conditions 
\[ \Delta_a(h_i)= \left\{ 
\begin{array}{cc}
1, & {\rm if} \; \; \; i=a, \\
0, & {\rm otherwise,} 
\end{array} \right. \] 
\[ \Delta_f(h_i)= \left\{ 
\begin{array}{cc}
1, & {\rm if} \; \; \; i=f, \\
0, & {\rm otherwise,} 
\end{array} \right. \] 
and the twisted Leibniz rule are well-defined on 
$\Q \langle {\cal YB}(G_2) \rangle .$ 
Let ${\cal Q}_a:=h_a^{-1} \circ \Delta_a$ and 
${\cal Q}_f:=h_f^{-1} \circ \Delta_f.$ 
The action of the simple reflections $s_a$ and $s_f$ on the Laurent polynomial 
ring $\Q[X_1^{\pm 1},X_2^{\pm 1}]$ is given by 
\[ s_a(X_1)=X_1, \; \; s_a(X_2)=X_1 X_2^{-1}, \] 
\[ s_f(X_1)=X_2, \; \; s_f(X_2)=X_1. \] 
Define the operators $\tau_a^{G_2}$ and $\tau_f^{G_2}$ acting on 
$\Q[X_1^{\pm 1} ,X_2^{\pm 1}]$ by 
\[ (\tau_a^{G_2}F)(X_1,X_2):= X_1\frac{F(X_1,X_2)-(s_aF)(X_1,X_2)}{X_2^2-X_1}, \] 
\[ (\tau_f^{G_2}F)(X_1,X_2):= X_2\frac{F(X_1,X_2)-(s_fF)(X_1,X_2)}{X_1-X_2}. \] 
The arguments as in the previous section show the following. 
\begin{pr}
\[ {\cal Q}_aF(\Theta_1,\Theta_2)=(\tau_a^{G_2}F)(\Theta_1,\Theta_2), \; \; 
{\cal Q}_fF(\Theta_1,\Theta_2)=(\tau_f^{G_2}F)(\Theta_1,\Theta_2). \]
\end{pr}
\begin{pr} {\rm 
There exists a natural homomorphism from 
$\Q \langle {\cal YB}(G_2) \rangle $ 
to the Nichols algebra 
${\cal B}_{G_2}$ obtained by $h_{\alpha}\mapsto 1+[\alpha],$ $\alpha\in \Psi_+.$ 
In other words, the $G_2$ Yang-Baxter relation holds in ${\cal B}_{G_2}.$} 
\end{pr} 
{\it Proof.} The Yang-Baxter relations give a set of 
relations among $[a],\ldots,[f]$ up to degree six. It is easy to check the 
compatibility of the quadratic relations and those from subsystems of type 
$A_2.$ The rest of cubic relations and the ones of higher degree can be verified 
by direct computation with help of the factorization of the 
braided symmetrizer, 
\cite{Ba}. ~\rule{3mm}{3mm} \medskip \\ 
The independent $W(G_2)$-invariant Laurent polynomials are given by 
\[ \phi_1(X_1,X_2)= X_1+X_1^{-1}+X_2+X_2^{-1}+X_1X_2^{-1}+X_1^{-1}X_2, \] 
\[ \phi_2(X_1,X_2)= X_1X_2+X_1^{-1}X_2^{-1}+X_1^2X_2^{-1}+X_1^{-1}X_2^2 + 
X_1^{-2}X_2 + X_1X_2^{-2} . \] 
The propositions above imply: 
\begin{theorem} We have $\phi_1(\Theta_1,\Theta_2)=
\phi_2(\Theta_1,\Theta_2)=6$ in the Nichols algebra ${\cal B}_{G_2},$ 
so the subalgebra of ${\cal B}_{G_2}$ generated by the images of 
the RSM-elements $\Theta_1^{G_2}$ and $\Theta_2^{G_2}$ is isomorphic to 
the Grothendieck ring of the flag variety of type $G_2.$ 
\end{theorem} 
\begin{de} Define the algebra ${\cal BE}(G_2)$ as an associative algebra
over $\Q$ with generators $\{a,b,c,d,e,f \}$ subject to the relations

$\bullet$ $($Commutativity$)$ $ad=da,~~be=eb,~~cf=fc;$

$\bullet$ $($Quadratic relations$)$ $ae=ec+ca,~ea=ce+ac,~fb=df+bd,~bf=fd+db;$
$$af=ba+cb+dc+ed+fe,~~fa=ab+bc+cd+de+ef;$$

$\bullet$ $($Quartic relations$)$ 
$$abac+acab+acbc=baca+cbca+caba,~
dfef+dedf+efdf=fded+fdfe+fefd,$$
$$abde+bcde+bcef+ecdb=cdbc+cdcd+decd+fdca,$$
$$bdce+edcb+edba+fecb=cbdc+dcdc+dced+acdf,$$
$\bullet$ $(G_2$ Yang--Baxter relation$)$ $abcdef=fedcba.$
\end{de}
\begin{con} The relations $\phi_1(\Theta_1,\Theta_2)=
\phi_2(\Theta_1,\Theta_2)=6$ are still valid in the algebra ${\cal BE}(G_2).$
\end{con}
\begin{rem} { \rm One can show that there exists the natural epimorphism of
algebras ${\cal BE}(G_2) \longrightarrow {\cal B}_{G_2},$  which has a
non-trivial kernel, however.
}
\end{rem}
\setcounter{section}{1}
\renewcommand{\thesection}{\Alph{section}}
\section*{Appendix}
\setcounter{theorem}{0}
\begin{de}
Let ${\cal YB}(B_n)$ be a group generated by the elements $\{ h_{ij}, 
g_{ij} \mid 1 \le i \not= j \le n \}$ and $ \{h_i \mid 1 \le i \le n \},$
subject to the following set of relations: \smallskip 

$\bullet$ $g_{ij}=g_{ji},$ $h_{ij}=h_{ji}^{-1},$ 

$\bullet$ $h_{ij}~h_{kl}=h_{kl}~h_{ij},~~g_{ij}~g_{kl}=g_{kl}~g_{ij},$~~ 
$h_{k}~h_{ij}=h_{ij}~h_{k},$ ~~$h_{k}~g_{ij}=g_{ij}~h_{k},$ 

if all $i,j,k,l$ are distinct; \smallskip 

$h_i~h_j=h_j~h_i,$ if $1 \le i, j \le n$;
~~~$h_{ij}~g_{ij}=g_{ij}~h_{ij},$  if $1 \le i < j \le n;$ \smallskip 

$\bullet$ {\rm ($A_2$ Yang-Baxter relations)}
           $$ ({\bf I})~~~h_{ij}~h_{ik}~h_{jk}=h_{jk}~h_{ik}~h_{ij},$$
           $$ ({\bf II})~~~h_{ij}~g_{ik}~g_{jk}=g_{jk}~g_{ik}~h_{ij},$$
           $$ ({\bf III})~~~h_{ik}~g_{ij}~g_{jk}=g_{jk}~g_{ij}~h_{ik},$$
           $$ ({\bf IV})~~~h_{jk}~g_{ij}~g_{ik}=g_{ik}~g_{ij}~h_{jk},$$
if $1\le i < j <k \le n;$ \smallskip 

$\bullet$  {\rm ($B_2$ quantum Yang-Baxter relation)}
$$h_{ij}~h_i~g_{ij}~h_j=h_j~g_{ij}~h_i~h_{ij},$$
if $ 1 \le i < j \le n.$
\end{de}
\begin{de}
Define the following elements in the group ${\cal YB}(B_n):$
\[ \Theta_j=(\prod_{i=j-1}^{1}h_{ij}^{-1})~h_j~(\prod_{i=1,i \ne j}^{n}g_{ij})
~h_j~(\prod_{k=n}^{j+1} h_{jk}), \]
for $1 \le j \le n. $ 
\end{de}
Proof of Theorem 2.3 (key lemma). 
We have to prove that $\Theta_i~\Theta_j~\Theta_i^{-1}=\Theta_j.$ It is enough
to consider the case $i < j.$

To begin with, it is convenient to introduce a bit of notation:

$$ A_i=\prod_{a=i-1}^{1}h_{ai}^{-1},~~~B_{i}=
\prod_{\scriptstyle a=1 \atop\scriptstyle a \not= i}^{n}g_{ai},~~~C_i=\prod_{a=n}^{i+1}h_{ia},$$
and in a similar way, we define $A_j,B_j$ and $C_j;$  

$$A_j^{'}=\prod_{\scriptstyle c=j-1 \atop\scriptstyle c \not= i}^{1}
h_{cj}^{-1},~
B_i^{'}=\prod_{\scriptstyle a=1 \atop\scriptstyle a \not= i,j}^{n}g_{ai},~~
B_j^{'}=\prod_{\scriptstyle c=1 \atop\scriptstyle c \not= i,j}^{n}g_{cj},~~
C_i^{'}=\prod_{\scriptstyle a=n \atop\scriptstyle a \not= j}^{i+1}h_{ia}.$$

Using this notation we can write
\[ \Theta_i~\Theta_j~\Theta_i^{-1}=A_{i}~h_{i}~B_{i}~h_{i}~C_{i}~A_{j}~h_{j}~
B_{j}~h_{j}~C_{j}~C_{i}^{-1}~h_{i}^{-1}~B_{i}^{-1}~h_{i}^{-1}~A_{i}^{-1}. \] 
The Lemma below describes commutation relations between elements we have 
introduced.
\begin{lem}
$(1a)$ ~~$C_{i}~A_{j}=A_{j}^{'}~C_{i}^{'};$

$(2a)$ ~~$C_{j}~C_{i}^{-1}=(C_{i}^{'})^{-1}~h_{ij}~C_{j};$

$(3a)$ ~~$C_{i}^{'}~B_{j}=B_{j}^{'}~g_{ij}~C_{i}^{'};$

$(4a)$ ~~$B_{i}~A_{j}^{'}=A_{j}^{'}~g_{ij}~B_{i}^{'};$

$(5a)$ ~~$C_{j}~B_{i}^{-1}=g_{ij}^{-1}~(B_{i}^{'})^{-1}~C_{j};$

$(6a)$ ~~$B_{j}^{'}~h_{ij}^{-1}~(B_{i}^{'})^{-1}=
B_{i}^{'}~h_{ij}^{-1}~(B_{j}^{'})^{-1};$

$(7a)$ ~~$A_{i}~A_{j}^{'}~h_{ij}^{-1}=A_{j}~A_{i};$

$(8a)$ ~~$g_{ij}~B_{j}^{'}~A_{i}^{-1}=A_{i}^{-1}~B_{j};$

$(9a)$ ~~$[h_i,B_j]=0,~[h_j,B_i]=0,~[h_j,A_i]=0,~[h_i,C_j]=0,$

$[h_i,A_{j}^{'}]=0,~[h_j,~C_{i}^{'}]=0,~[A_i,C_j]=0.$

\end{lem} 

Using relations $(1a)$ and $(2a),$ and commutativity relations $(9a)$, we can 
write
\[ \Theta_i~\Theta_j~\Theta_{i}^{-1}=A_{i}~h_{i}~{\bf B_{i}~A_{j}^{'}}~h_{i}
~h_{j}~{\bf C_{i}^{'}~B_{j}~(C_{i}^{'})^{-1}}~h_{j}~h_{ij}^{-1}~h_{i}^{-1}
~{\bf C_{j}~B_{i}^{-1}}~h_{i}^{-1}~A_{i}^{-1}. \] 
Now we are going to apply the relations $(3a),(4a)$ and $(5a)$ respectivly 
to the market terms to reduce the above expression to the following form
\[ A_{i}~A_{j}^{'}~h_{i}~g_{ij}~B_{i}~h_{i}~h_{j}~B_{j}^{'}
~{\bf g_{ij}~h_{j}~h_{ij}^{-1}~h_{i}^{-1}}~g_{ij}^{-1}~(B_{i}{'})^{-1}
~C_{j}~h_{i}^{-1}~A_{i}^{-1}. \] 
To the market terms we can apply the $B_2$-Yang-Baxter relation presented in 
an equivalent form  $g_{ij}~h~{j}~h_{ij}^{-1}~h_{i}^{-1}=
h_{i}^{-1}~h_{ij}^{-1}~h~{j}~g_{ij},$  and after that do cancellations 
of $h_i$ and $g_{ij}.$  As a result we will have
\[ \Theta_i~\Theta_j~\Theta_i^{-1}=
A_{i}~A_{j}^{'}~h_{i}~g_{ij}~B_{i}^{'}~h_{j}~
{\bf B_{j}^{'}~h_{ij}^{-1}(B_{i}^{'})^{-1}}~h_{j}~
C_{j}~h_{i}^{-1}~A_{i}^{-1}. \] 
The next step is to apply to the bold terms the relation $(6a),$ and rewrite 
the above expression in the following form:
\[ A_{i}~A_{j}^{'}~{\bf h_{i}~g_{ij}~h_{j}~h_{ij}^{-1}}~B_{j}^{'}~h_{j}~C_{j}~
h_{i}^{-1}~A_{i}^{-1}. \] 
Now we can apply to the market terms the $B_2$-Yang-Baxter relation again, 
but this time written in the form $h_{i}~g_{ij}~h~{j}~h_{ij}^{-1}=
h_{ij}^{-1}~h~{j}~g_{ij}~h_{i},$ and  after the cancellation of $h_{i},$ to 
obtain
\[ {\bf A_{i}~A_{j}^{'}~h_{ij}^{-1}}~g_{ij}~B_{j}^{'}~h_{j}~C_{j}~A_{i}^{-1}. \]
Now applying the relation $(7a)$ to the market terms, and using the fact that 
$C_j$ and $A_i$ commute, see relations $(9a),$  we can write
\[ \Theta_i~\Theta_j~\Theta_i^{-1}=A_{j}~A_{i}~h_{j}
~{\bf g_{ij}~B_{j}^{'}~A_{i}^{-1}}~h_{j}~C_{j}. \] 
Finally, applying the relation $(8a)$ to the market tems, after cancellations 
we will have
\[ \Theta_i~\Theta_j~\Theta_i^{-1}=A_{j}~h_{j}~B_{j}~h_{j}~C_{j}=\Theta_{j}. 
\; \; \; \blacksquare \] 
It remains to prove the commutativity relations listed in Lemma A.3.
\bigskip \\
$\bullet$ {\it Proof of $(1a)$}:

Note that if $j=i+1,$ then $C_i~A_j=1=A_{j}^{'}~C_i^{'}.$ So we will assume 
that $j-i \ge 2.$ Under this assumption, we can write 
\[ C_{i}~A_{j}=h_{i,n}\cdots {\bf h_{i,j}~h_{i,j-1}} \cdots h_{1,i}~
{\bf h_{j-1,j}^{-1}} \cdots h_{ij}^{-1} \cdots h_{1,j}^{-1}. \] 
Using local commutativity relations, see Definition~3.1, 
$(2),$ we can move the factor 
${\bf h_{j-1,j}^{-1}}$ to the left until we have touched on the factor 
${\bf h_{i,j-1}}.$ As a result, we will come up with the triple product
${\bf h_{i,j}~h_{i,j-1}~ h_{j-1,j}^{-1}}.$ 
Now we can apply the $A_2$-Yang-Baxter relation, see Definition A.1, 
$({\bf I}),$ in an equivalent form $h_{jk}^{-1}~h_{ij}~h_{ik}=
h_{ik}~h_{ij}~h_{jk}^{-1},$ 
and move the factor ${\bf h_{j-1,j}}$ to the left most position, to obtain
\[ C_{i}~A_{j}=h_{j-1,j}^{-1}~h_{i,n} \cdots h_{i,j-1}~{\bf h_{ij}~h_{i,j-2}} 
\cdots h_{i,i+1}~{\bf h_{j-2,j}^{-1}} \cdots h_{ij}^{-1} \cdots h_{1,j}^{-1}. \] 
In a similar fashion as above, we can move the factor ${\bf h_{j-2,j}^{-1}}$ 
to the left until we have touched on the factor ${\bf h_{i,j-2}}.$ Now we can 
apply the $A_2$-Yang-Baxter relation $({\bf I})$ in the form presented above, 
to the triple product $ {\bf h_{ij}~h_{i,j-2}~h_{j-2,j}^{-1}},$ and move to 
the left the factor ${\bf h_{j-2,j}^{-1}}.$ We can continue this procedure 
until the factor  ${\bf h_{ij}}$ will touch the factor $h_{ij}^{-1}$. 
After cancellation and moving to the left the product 
$h_{i-1,j}^{-1} \cdots h_{1,j}^{-1},$ we will come to
the product  $A_{j}^{'}~C_i^{'}.$ $\blacksquare$
\bigskip \\ 
$\bullet$ {\it Proof of $(2a)$} is similar to that of $(1a).$

By definition, $C_{j}~C_{i}^{-1}= h_{j,n} 
\cdots {\bf h_{j,j+1}}~h_{i,i+1}^{-1} 
\cdots {\bf h_{ij}^{-1}~h_{i,j+1}^{-1}} \cdots h_{i,n}^{-1}.$ Using local 
commutativity relations, see Definition~3.1, $(2),$ we can move the factor 
${\bf h_{j,j+1}}$ to the right until we have touched on the factor 
${\bf h_{ij}^{-1}}.$ As a result, we will come up with the triple product
${\bf h_{j,j+1}~h_{ij}^{-1}~ h_{i,j+1}^{-1}}.$ 

Now we can apply the $A_2$-Yang-Baxter relation $({\bf I}),$ see 
Definition A.1,  in an equivalent form $h_{jk}~h_{ij}^{-1}~h_{ik}^{-1}=
h_{ik}^{-1}~h_{ij}^{-1}~h_{jk},$ and move the 
factor  ${\bf h_{j,j+1}}$ to the right most position and 
$h_{i,j+1}^{-1}$ to the left most position, to obtain 
$$C_{j}~C_{i}^{-1}=h_{i,j+1}^{-1}~h_{j,n} \cdots {\bf h_{j,j+2}}~
h_{i,i+1}^{-1} \cdots {\bf h_{ij}^{-1}~h_{i,j+2}^{-1}} \cdots h_{i,n}^{-1}~
h_{j,j+1}.$$
Now we can move the factor ${\bf h_{j,j+2}}$ to the right until we have 
touched on the factor ${\bf h_{ij}^{-1}},$ and apply the $A_2$-Yang-Baxter 
relation $({\bf I})$ in the form mentioned above, to the triple product
$${\bf h_{j,j+2}~h_{ij}^{-1}~h_{i,j+2}^{-1}}.$$
Just as before, we will come to the equality
$$C_{j}~C_{i}^{-1}=h_{i,j+1}^{-1}~h_{i,j+2}^{-1}~h_{j,n} \cdots 
{\bf h_{j,j+3}}~h_{i,i+1}^{-1} \cdots {\bf h_{ij}^{-1}~h_{i,j+3}^{-1}} 
\cdots h_{i,n}^{-1}~h_{j,j+2}~h_{j,j+1}.$$
Repeating this procedure we will come to 
$$\prod_{\scriptstyle a=i+1 \atop\scriptstyle a \not= j}^{n}h_{ia}^{-1}
~h_{ij}^{-1}~\prod_{c=n}^{j+1}h_{jc}.$$ 
$\bullet$ {\it Proof of $(3a)$} is similar to that of $(1a)$ and $(2a),$ but
this time we have to use $A_2$-Yang-Baxter relation $({\bf II}).$ 

By definition, $C_{i}^{'}~B_{j}=h_{i,n} \cdots {\widehat h_{ij}} \cdots 
{\bf h_{i,i+1}}~g_{1,j} \cdots {\bf g_{ij}~g_{i+1,j}} \cdots g_{n,j}.$
We can move the factor ${\bf h_{i,i+1}}$ to the right until we have touched 
on the factor ${\bf g_{ij}}.$  Now we can apply the $A_2$-Yang-Baxter relation 
${\bf II}$ to the triple product ${\bf h_{i,i+1}~ g_{ij}~g_{i+1,j}}.$ As a 
result, we will come to an equality
$$C_{i}^{'}~B_{j}=h_{i,n} \cdots {\widehat h_{ij}} \cdots 
{\bf h_{i,i+2}}~g_{1,j} \cdots g_{i-1,j}~g_{i+1,j}~{\bf g_{ij}~g_{i+2,j}} \cdots g_{nj}~h_{i,i+1}.$$
Now we can again move the factor ${\bf h_{i,i+2}}$ to the right until we have
touched on the factor ${\bf g_{ij}},$ and then apply the $A_2$-Yang-Baxter 
relation ${\bf II}$ to the triple product 
${\bf h_{i,i+2}~ g_{ij}~g_{i+2,j}}.$ The result can be written as follows
$$C_{i}^{'}~B_{j}=h_{i,n} \cdots {\widehat h_{ij}} \cdots 
{\bf h_{i,i+3}}~g_{1,j} \cdots g_{i-1,j}~g_{i+1,j}~g_{i+2,j}~
{\bf g_{ij}~g_{i+3,j}} \cdots g_{n,j}~h_{i,i+2}~h_{i,i+1}.$$
Repeating this procedure we can move the all $h_{i,k}, ~k \not=j,$ 
to the right through the factor 
$g_{ij}.$ $\blacksquare$
\bigskip \\ 
$\bullet$ {\it Proof of $(4a)$} is a very similar to that of $(3a),$ but
this time we have to use $A_2$-Yang-Baxter relation $({\bf IV})$ in the form 
$$g_{ik}~g_{ij}~h_{jk}^{-1}=h_{jk}^{-1}~g_{ij}~g_{ik}.$$ 
By definition 
$$B_{i}~A_{j}^{'}=g_{1,i} \cdots {\bf g_{i,j-1}~g_{ij}} \cdots 
g_{i,n}~{\bf h_{j-1,j}^{-1}} \cdots {\widehat h_{ij}^{-1}} 
\cdots h_{1,j}^{-1}.$$
We can move the factor ${\bf h_{j-1,j}^{-1}}$ to the left until we have touched
on the factor ${\bf  g_{ij}}.$ Then we can apply the $A_2$-Yang-Baxter 
relation $({\bf IV})$ in the form presented above, and transform the result 
to the following form
$$h_{j-1,j}^{-1}~g_{1,i} \cdots {\bf g_{i,j-2}~g_{ij}} \cdots 
g_{i,n}~{\bf h_{j-2,j}^{-1}} \cdots {\widehat h_{ij}^{-1}} 
\cdots h_{1,j}^{-1}~g_{i,j-1}.$$
Now we can move the factor ${\bf h_{j-2,j}^{-1}}$ to the left until we have 
touched on the factor ${\bf g_{ij}},$ and apply the Yang-Baxter relation 
$({\bf IV})$ to the triple product ${\bf g_{i,j-2}~g_{ij}~h_{j-2,j}^{-1}},$ 
and so on. As a final result we will come to the RHS of the equality $(4a).$ 
$\blacksquare$
\bigskip \\ 
$\bullet$ {\it Proof of $(5a)$} is a very similar to that of $(4a),$ but
this time we have to use $A_2$-Yang-Baxter relation $({\bf IV})$ in the 
following form 
$$h_{jk}~g_{ik}^{-1}~g_{ij}^{-1}=g_{ij}^{-1}~g_{ik}^{-1}~h_{jk}.$$
We will give only an outlook of the proof and leave details to the reader.\\
By definition 
$$ C_{j}~B_{i}^{-1}=h_{jn} \cdots {\bf h_{j,j+1}}~g_{in}^{-1} \cdots 
{\bf g_{i,j+1}~g_{ij}^{-1}} \cdots g_{1,j}^{-1}.$$
Therefore we can apply to the market terms the Yang-Baxter relation 
$({\bf IV})$  in the form presented above and write
$$ C_{j}~B_{i}^{-1}=h_{jn} \cdots { \bf h_{j,j+2}}~g_{in}^{-1} \cdots 
{\bf g_{i,j+2}^{-1}~g_{ij}^{-1}}~g_{i,j+1}^{-1} \cdots g_{1,j}^{-1}~h_{j,j+1}.$$
We can continue by applying the Yang-Baxter relation $({\bf IV})$  to the 
market factors, and so on. As a final result we obtain the RHS of the 
equality $(5a).$ $\blacksquare$
\bigskip \\
$\bullet$ {\it Proof of $(6a)$} runs in the same way as that of $(5a),$ but
this time we have to use $A_2$-Yang-Baxter relation $({\bf II})$ in the 
following form 
$$g_{jk}~h_{ij}^{-1}~g_{ik}^{-1}=g_{ik}^{-1}~h_{ij}^{-1}~g_{jk}.$$
Again we will give only an outlook on the proof and leave details to the 
reader.\\
By definition
$$B_{j}^{'}~h_{ij}^{-1}~(B_{i}^{'})^{-1}=g_{1,j} \cdots {\hat {g}_{ij}} 
\cdots {\bf g_{j,n}~h_{ij}^{-1}~g_{i,n}^{-1}} \cdots {\hat {g}_{ij}} \cdots 
g_{1,i}^{-1}.$$
Now we can apply the Yang-Baxter relation $({\bf II})$ in the form presented 
above to the market terms to conclude that
$$ B_{j}^{'}~h_{ij}^{-1}~(B_{i}^{'})^{-1}=g_{i,n}^{-1}~g_{1,j} \cdots 
{\hat {g}_{ij}} \cdots {\bf g_{j,n-1}~h_{ij}^{-1}~g_{i,n-1}^{-1}} \cdots 
{\hat {g}_{ij}} \cdots g_{1,i}^{-1}~g_{j,n}.$$
Now we can continue and apply the Yang-Baxter relation $({\bf II})$ to the 
market terms, and so on. As a final step we have come to the RHS of the 
equality $(6a).$ $\blacksquare$
\bigskip \\ 
$\bullet$ {\it Proof of $(7a)$} runs in the same way as that of the previous 
one, but this time we have to use $A_2$-Yang-Baxter relation $({\bf I})$ in 
the  following form 
$$h_{ij}^{-1}~h_{ik}^{-1}~h_{jk}^{-1}=h_{jk}^{-1}~h_{ik}^{-1}~h_{ij}^{-1}.$$
Again we will give only an outlook on the proof and leave details to the 
reader.\\
By definition
$$A_i~A_j~h_{ij}^{-1}=h_{i-1,i}^{-1} \cdots {\bf h_{1,i}^{-1}}h_{j-1,j}^{-1}
\cdots {\widehat h}_{ij}^{-1} \cdots {\bf h_{1,j}^{-1}~h_{ij}^{-1}}.$$
Therefore we can apply to the market terms the Yang-Baxter relation 
$({\bf I})$ in the form presented above and write
$$A_i~A_j~h_{ij}^{-1}=h_{i-1,i}^{-1} \cdots {\bf h_{2,i}^{-1}}h_{j-1,j}^{-1}
\cdots {\widehat h}_{ij}^{-1} \cdots {\bf h_{2,j}^{-1}~h_{ij}^{-1}}~h_{1,j}^{-1}~h_{1,i}^{-1}.$$
Now we can continue and apply the Yang-Baxter relation $({\bf I})$ to the 
market terms, and so on. As a final step we have come to the RHS of the 
equality $(7a).$ $\blacksquare$
\bigskip \\ 
$\bullet$ {\it Proof of $(8a)$} runs in the same way as that of the previous 
one, but this time we have to use $A_2$-Yang-Baxter relation $({\bf II}).$ 
We will give only an outlook on the proof and leave details to the reader.\\
By definition
$$A_{i}^{-1}~B_{j}=h_{1,i} \cdots {\bf h_{i-1,i}}~g_{1,j} \cdots 
{\bf g_{i-1,j}~g_{ij}} \cdots g_{nj}.$$
Therefore we can apply to the market terms the Yang-Baxter relation 
$({\bf II})$ in the form presented above and write
$$A_{i}^{-1}~B_{j}=h_{1,i} \cdots {\bf h_{i-2,i}}~g_{1,j} \cdots 
{\bf g_{i-2,j}~g_{ij}}~g_{i-1,j}~ \cdots g_{nj}~h_{i-1,i}.$$
Now we can continue and apply the Yang-Baxter relation $({\bf II})$ to the 
market terms, and so on. As a final step we have come to the RHS of the 
equality $(8a).$
\bigskip \\ 
$${\bf Acknowledgements}$$ 
The authors would like to thank Yuri Bazlov for fruitful discussions and help
with computation of the Hilbert series of the Nichols-Woronowicz algebras
${\cal B}_{B_3}$ and that ${\cal B}_{G_2}.$ 
Both of the authors were 
supported by Grant-in-Aid for Scientific Research.

Research Institute for Mathematical Sciences, \\
Kyoto University, \\ 
Sakyo-ku, Kyoto 606-8502, Japan \\ 
e-mail: kirillov@kurims.kyoto-u.ac.jp 
\medskip \\
Department of Mathematics, \\
Kyoto University, \\ 
Sakyo-ku, Kyoto 606-8502, Japan \\ 
e-mail: maeno@math.kyoto-u.ac.jp


\begin{thebibliography}{11}

\bibitem{Ba} Bazlov Y., {\it Nichols-Woronowicz algebra model for Schubert 
calculus on Coxeter groups,} preprint, math.QA/0409206. 

\bibitem{BGG} Bernstein I.N, Gelfand I.M. and Gelfand S.I., {\it Schubert 
cells and cohomology of the space $G/P$}, Russian Math. Survey {\bf 28} (1973),
1-26.

\bibitem{FK} Fomin S.V. and Kirillov A.N., {\it Quadratic algebras, Dunkl
elements, and Schubert calculus}, Advances in Geometry, 147-182, 
Progress in Math. {\bf 172}, Birkhauser, Boston, 1998. 

\bibitem{GL} Givental A. and Lee Y.-P., 
{\it Quantum $K$-theory on flag manifolds, finite-difference Toda lattices and 
quantum groups,} Invent. Math., {\bf 151} (2003), 193-219. 


\bibitem{Kir2} Kirillov A.N., {\it On some quadratic algebras: Jucys-Murphy 
and Dunkl elements}, Calogero-Moser-Sutherland models, (Montr\'eal, PQ, 
March 10-15, 1997), CRM Series in Mathematical Physics, Springer, New York, 
2000, 231-248. 

\bibitem{Kir3} Kirillov A.N., {\it Quantum Grothendieck polynomials,}
   Algebraic Methods and q-Special Functions (Montr\'eal, PQ, May 13-17, 1996), 
   CRM Proc. Lecture Notes {\bf 22}, AMS, 1999, 215-226;
   q-alg/9610034. 

\bibitem{KM} Kirillov A.N. and Maeno T., {\it Noncommutative algebras related 
with Schubert calculus on Coxeter groups}, Europian Journ. of Combin.
{\bf 25}, (2004), 1301-1325.


\bibitem{KK1} Kostant B. and Kumar S., {\it The nil Hecke ring and cohomology 
of $G/$ for a Kac-Moody group $G$ }, Advances in Math. {\bf 62} (1986), 
187-237.

\bibitem{KK2} Kostant B. and Kumar S., {\it T-equivariant K-theory of 
generalized flag varieties}, Journ. Differential Geometry {\bf 32} (1990),
549-603.   

\bibitem{L1} Lascoux A., {\it Anneau de Grothendieck de la vari\'et\'e de 
drapeaux}, The Grothendieck Festschrift, Vol. III, 1--34, Progress in Math., 
{\bf 88}, Birkh\"auser Boston, Boston, MA, 1990. 


\bibitem{LS1} Lascoux A. and Sch\"utzenberger M.-P., {\it Symmetry and flag 
manifolds},  Invariant theory (Montecatini, 1982), 118--144,
Lecture Notes in Math., {\bf 996}, Springer, Berlin, 1983.

\bibitem{LS2} Lascoux A. and Sch\"utzenberger M.-P., {\it Polynomes 
d\'e Schubert},  C. R. Acad. Sci. Paris  S\'er. I Math. {\bf 294} (1982), 
{\bf  13}, 447--450. 

\bibitem{len1} Lenart C., {\it The K-theory of the flag variety and the 
Fomin-Kirillov quadratic algebra},  J. Algebra {\bf 285} (2005), 120-135.

\bibitem{len2} Lenart C. and Sottile F., {\it A Pieri-type formula for the 
K-theory of a flag manifold}, preprint, math.CO/0407412. 
   
\bibitem{LY} Lenart C. and Yong A., {\it Lecture notes on the K-theory 
of the flag variety and the Fomin-Kirillov quadratic algebra,} 
\verb|http://www.math.umn.edu/~ayong/|


\bibitem{Po} Postnikov, A., {\it On a quantum version of Pieri's formula,} 
Advances in Geometry, 371-383, 
Progress in Math. {\bf 172}, Birkhauser, Boston, 1998. 

\end{thebibliography}
\end{document}